\documentclass[preprint,review,12pt,3p]{elsarticle}
\usepackage{amsmath}
\usepackage{graphicx,amssymb,psfrag,subfigure,color}
\usepackage{enumitem}
\usepackage{epsfig}
\usepackage{bm}
\usepackage{cases}
\usepackage{upgreek}
\usepackage{array}
\usepackage{lineno}
\usepackage[colorlinks]{hyperref}
\usepackage{epstopdf}
\usepackage{caption}
\usepackage{booktabs}
\usepackage{natbib}
\usepackage[table]{xcolor}
\usepackage[section]{placeins}

\makeatletter
\renewcommand{\maketag@@@}[1]{\hbox{\m@th\normalsize\normalfont#1}}%
\makeatother
\newcommand{\beq}{\begin{equation}}
\newcommand{\eeq}{\end{equation}}
\newcommand{\bpm}{\begin{pmatrix}}
\newcommand{\epm}{\end{pmatrix}}
\newcommand{\beqa}{\begin{eqnarray}}
\newcommand{\eeqa}{\end{eqnarray}}
\newcommand{\beqas}{\begin{eqnarray*}}
\newcommand{\eeqas}{\end{eqnarray*}}

\newcommand{\pdhfrac}[2]{\mathchoice{\frac{#1}{#2}}{#1/#2}{#1/#2}{#1/#2}}

\newcommand{\pd}[2]{\pdhfrac{{\partial}#1}{{\partial}#2}}

\newcommand{\vD}{\textbf{D}}

\newcommand{\vx}{\textbf{x}}

\newcommand{\vy}{\textbf{y}}
\newcommand{\vY}{\textbf{Y}}

\newcommand{\vt}{\textbf{t}}

\newcommand{\vn}{\textbf{n}}

\newcommand{\vu}{\textbf{u}}
\newcommand{\vf}{\textbf{f}}
\newcommand{\vF}{\textbf{F}}

\newcommand{\vJ}{\textbf{J}}

\newcommand{\eps}{\epsilon}

\newcommand{\la}{\lambda}

\newcommand*{\dif}{\mathop{}\!\mathrm{d}}

\def\XXint#1#2#3{{\setbox0=\hbox{$#1{#2#3}{\int}$ }
\vcenter{\hbox{$#2#3$ }}\kern-.6\wd0}}

\biboptions{sort&compress}

\hypersetup{
linkcolor=cyan,
anchorcolor=black,
filecolor=black,
urlcolor=black,
citecolor=cyan
}
\begin{document}

\begin{frontmatter}

\title{Optimisation of spatially varying orthotropic porous structures based on conformal mapping}


\author[DUT1,DUT2]{Shaoshuai Li}
\author[DUT1,DUT2,DUT3]{Yichao Zhu\corref{cor1}}
\ead{yichaozhu@dlut.edu.cn}
\author[DUT1,DUT2,DUT3]{Xu Guo\corref{cor1}}
\ead{guoxu@dlut.edu.cn}
\cortext[cor1]{Corresponding author}
\address[DUT1]{State Key Laboratory of Structural Analysis for Industrial Equipment, Department of Engineering Mechanics, Dalian University of Technology, Dalian, 116023, P. R. China}
\address[DUT2]{International Research Center for Computational Mechanics, Dalian University of Technology}
\address[DUT3]{Ningbo Institute of Dalian University of Technology, No.26 Yucai Road, Jiangbei District, Ningbo, 315016, P. R. China}

\begin{abstract}
  In this article, a compliance minimisation scheme for designing spatially varying orthotropic porous structures is proposed. With the utilisation of conformal mapping, the porous structures here can be generated by two controlling field variables, the (logarithm of) the local scaling factor and the rotational angle of the matrix cell, and they are interrelated through the Cauchy-Riemann equations. Thus the design variables are simply reduced to the logarithm values of the local scaling factor on selected boundary points. Other attractive features shown by the present method are summarised as follows. Firstly, with the condition of total differential automatically met by the two controlling field variables, the integrability problem which necessitates post-processing treatments in many other similar methods can be resolved naturally. Secondly, according to the maximum principle for harmonic functions, the minimum feature size can be explicitly monitored during optimisation. Thirdly, the rotational symmetry possessed by the matrix cell can be fully exploited in the context of conformal mapping, and the computational cost for solving the cell problems for the homogenised elasticity tensor is maximally abased. In particular, when the design domain takes a rectangle shape, analytical expressions for the controlling fields are available. The homogenised results are shown, both theoretically and numerically, to converge to the corresponding fine-scale results, and the effectiveness of the proposed work is further demonstrated with more numerical examples.
\end{abstract}

\end{frontmatter}


\section{Introduction}
Porous structures have demonstrated their exceptional performances in various engineering fields \cite{Lakes_Nature1993, Sigmund_JMPS1997, Lu_Science1999, Liu_JAM2017, Krauss_Pqe1999, Liu_Pra2015}, and nowadays the manufacturability of porous structures is increasingly enhanced, thanks to the rapid development of additive manufacturing technology. This naturally gives rise to the demand for design tools for porous configurations, and the major challenge stems from the fact that a porous structure is generally associated with more than one length scale, normally a macroscopic scale characterised by its overall size and a microscale identified by the size of its constituting microstructural members. In consideration of such a multiscale feature, analysis of the behaviour of porous structures seems to be constrained down to the microscopic level, which, albeit a number of effective works \cite{Alexandersen_CMME2015, Liu_JAM2017, Liu_CMME2020}, often poses high demands on computation resources.

To resolve the issue of efficiency, homogenisation treatments for porous configurations are also under researchers' attention. The aim is to decouple the original problem into two subproblems formulated on separate scales. The subproblem on the microscale, which is sometimes called the cell problem, is often localised of a macroscopic material point (Gauss point in numerical implementation) and designated to estimate the equivalent property. The macroscopic subproblem is coarse-grained, and solved with homogenised) properties fed from the solutions to the microscale subproblem. When the microscopic configuration of a porous material becomes spatially periodic, the solution based on the homogenised elasticity tensor sees its rigorousness via perturbation theory, i.e., the scale-separated solution asymptotically converges to that of the original multiscale problem as the length scale ratio tends to zero. This forms the backbone of the asymptotic-analysis-based homogenisation (AABH) method \cite{Bensoussan_Book1978, cioranescu_Book1999, Pavliotis_book2008}. Bends{\o}e and Kikuchi \cite{Bendsoe_CMME1988} pioneerly implemented the AABH method for the topology design of continuum structures. In their frameworks, any macroscopic point in the prescribed domain is assumed being composed of infinitesimal and periodically distributed microstructures. Structural optimisation is manifested through varying the microscopic parameters representing the geometry of the cell. By representing the microstructural configuration with a density field and assuming an explicit relation between the density and the assumed isotropic homogenised material property, the solid isotropic material with penalisation (SIMP) framework was introduced \cite{Bendsoe_so1989, Zhou_CMME1991}, then intensively developed and popularly used, because of its simplicity and effectiveness. Furthermore, Sigmund \cite{Sigmund_IJSS1994} proposed an inverse homogenisation approach to design the unit cell configurations with desirable elastic properties and carried out a series of microstructure design work, such as piezoelectric sensors \cite{Sigmund_JMR1998}, optimising thermal expansion coefficient \cite{Sigmund_JMPS1997}, and extremal elastic properties material \cite{Sigmund_2000JMPS} (more details in the monograph \cite{Bendsoe_book2003}). Then concurrent optimisation schemes \cite{Rodrigues_Smo2002, Coelho_SMO2008, LiuLing_CS2008} searching among both the macroscopic structural topologies and the microstructural configurations were proposed, and were then extended in several topics (e.g., \cite{Niu_SMO2009, Deng_SMO2013, Deng_SMO2017}, etc.). It is worth noting, however, that aforementioned works are devised based on the homogenisation results which are actually rigorous only for periodic configurations.

Compared with periodic porous structures, configurations infilled with spatially varying cells have greater design freedom. Early-stage attempts on the design and optimisation of graded porous structures have been made \cite{Zhou_JMS2008, Radman_JMS2012, Radman_CMS2014, Wang_CMME2017, Cheng_CMME2019}, but challenging issues persist. Firstly, only the interior layouts of the cells are permitted to vary, while the cell shape is kept fixed. Secondly, the microstructural variation must be in alignment with the directions that are orthogonal to one facet of the constituting cells. Thirdly, in most cases, spatial changes in cells are accomplished by (discretely) varying parameters controlling the cell structures. As a result, smooth connections between cells can not be ensured, as the parameter values are essentially piecewise constants.

Recently, further investigation into this subject has attracted wide attention, mainly based on conformal mapping \cite{Vogiatzis_CMME2018, Jiang_FME2019, Allaire_CMA2019, Xie_CAD2020} and the projection method \cite{Groen_IJNME2018, Groen_CMME2019, Groen_CMME2020, Perle_JCP2020} to generate configurations infilled with smoothly varying cells. The reason that conformal mapping is so attractive, as will be seen here, is that orthotropic porous configurations offer a natural choice for compliance optimisation \cite{Groen_IJNME2018}. Under these novel frameworks, the limiting issues outlined in the previous paragraph, such as limited cell orientation and smoothness, are properly addressed. Vogiatzis et al. \cite{Vogiatzis_CMME2018} considered using filling the unit cells with specific elastic properties to fill an arbitrary surface, but the issue of optimisation is not discussed. Jiang et al. \cite{Jiang_FME2019} discussed optimising porous structures with conformal mapping, but homogenised elasticity tensor is maintained to be isotropic. As a result, the rotational effect of the infilling cells is not considered. The projection-based method, which is also termed as a ``de-homogenisation'' method by some of its main contributors \cite{Groen_CMME2020}, sees its origin in the work by Pantz and Trabelsi \cite{Pantz_SIAM_JCO2008}. In this approach, the homogenisation-based topology optimisation is carried out on a coarse grid first. Then the optimised results are projected onto smooth graded and nearly orthotropic porous structures with a desired resolution, i.e. de-homogenisation. Here the projection method effectively offers a post-processing framework. The reason is briefed as follows. In order to fully resolve a porous configuration, de-homogenisation relies on the determination of a smooth mapping function, which indicates the spatial variation of the structure. But the design variables for compliance optimisation are the rotational angles of the matrix cells on pixels, which are related to the spatial derivatives of the mapping function. As a result, one has to integrate optimised results so as to finally describe the structure on the fine scale. However, with the conditions of total differential not met in general, the integrability of the results can not be ensured, and the projection method is then employed to ``project'' the derivatives onto a total differential form, so as to approximate the desired mapping function.

In parallel with the novel frameworks mentioned above, Zhu et al. \cite{Zhu_JMPS2019} proposed an AABH-plus-based optimisation scheme to design graded porous structures. The method consists of three modules: multiscale representation, AABH plus calculation, and compliance optimisation. For representation, a smooth and continuous function is introduced to map a graded porous structure to an artificial periodic structure. Then the moving morphable components/voids (MMC/MMV) approach \cite{Guo_JAM2014, Zhang_SMO2015, Zhang_CM2016,Zhang_JAM2017,Zhang_CMME2017} is adopted to describe the unit cell configuration. Thus the multiscale topology description function is formed through function composition \cite{Liu_JAM2017, Zhu_JMPS2019}. The AABH plus part is to (rigorously) derive the corresponding scale-separated form, by means of asymptotic analysis. For optimisation, the parameters controlling the (composite) topology description function (TDF) become the design variables. Upon modification, such treatments of representing graded microstructural configurations with composite TDFs actually provide a platform, on which existing descriptions of porous configurations get formulated in a unified manner \cite{Xue_CMME2020}. At its present stage, the AABH-based framework has yet realised its full potential due to the issues detailed as follows. Firstly, unlike periodic structures, cell configurations now differ point by point in the coarse-grained design domain. This means the microstructural cell problems should be computed as many times as the number of macroscopic finite elements, leading to a huge computational burden. This challenging issue is somehow mitigated by adopting a linearisation treatment \cite{Zhu_JMPS2019} or a zoning strategy \cite{Xue_Smo2020}. Secondly, the use of composite TDFs enhances the describability of graded porous configurations. Then how to manage such strong describability becomes a critical issue. For compliance optimisation, in particular, configurations infilled with orthotropic microstructures sometimes suffice in netting a nearly optimised solution for two-dimensional cases \cite{Groen_IJNME2018}. Therefore, the identification of effective porous configurations for optimisation is an issue worth further investigation.

To this end, a compliance optimisation framework among orthotropic porous configurations is proposed in this article. Here the macroscopic mapping function monitoring the gradual change of the matrix cell is restricted to holomorphic functions. Consequently, the porous configurations of interest are generated through gradually rescaling and rotating the matrix cell in space. As from the conformal mapping rule, the logarithm of the scaling factor denoted by $\ln\lambda$ here and the rotation angle denoted by $\theta$ here should just correspond to the real and imaginary parts of a holomorphic function. This greatly benefits the realisation of the resulting compliance optimisation, summarised as follows. Firstly, the porous configurations generated through holomorphic mapping functions maintain smooth connectivity naturally, and are orthotropic provided that the matrix cells are orthotropic. Secondly, the equivalent elasticity tensors should be computed from the cell problems only for the cases of the matrix cell, and the actual onsite (homogenised) elasticity tensors can be obtained by rotating the corresponding matrix cells implied by the angle $\theta$. Hence the homogenised results here are expected (and are also shown) to be asymptotically consistent with that from fine-scale calculations. But the computational efficiency should be enhanced greatly, compared to early-stage AABH-plus-based schemes \cite{Zhu_JMPS2019, Xue_Smo2020, Xue_CMME2020}. It is shown that the optimisation for two-dimensional cases can be brought down roughly below 100 seconds on a desktop computer. Thirdly, the two design controllers $\ln\lambda$ and $\theta$ are interrelated through Cauchy-Riemamm equations. The number of the design controller, which is chosen to be $\lambda$ here, is then reduced to 1. Fourthly, as $\ln\lambda$ must satisfy two-dimensional Laplace's equation, $\ln\lambda$ is fully determined with its prescribed boundary values. Thus the actual design variables are $\ln\lambda$ evaluated at a collection of boundary points. Compared to the SIMP-based schemes, the number of design variables is much fewer. Fifthly, the design controller $\ln\lambda$ is a harmonic function, whose behaviour is guided by the maximum principle, i.e. the extremal values of $\ln\lambda$ must be attained on the domain boundaries. Thus the minimum size of microstructural members can be explicitly controlled through the design variables, the boundary values of $\ln\lambda$. This is especially important when there are requirements on the minimal printable sizes in additive manufacturing. Finally, the rotational angle $\theta$ here is automatically continuous in space. This is in contrast with the projection-based approaches \cite{Groen_IJNME2018, Allaire_CMA2019, Xie_CAD2020}, where $\cos\theta$ or $\sin\theta$ serves as the actual design parameters, and a jump of angle $\pi$ between neighbouring pixels can not be fully avoided without imposing extra conditions.

To the present stage, features that distinguish the present work from others in the frontier, mainly the projection-based methods \cite[e.g.,][]{Groen_IJNME2018, Groen_CMME2019, Groen_CMME2020, Allaire_CMA2019, Perle_JCP2020}, are specified in brief, while further elaboration over this point is conducted in the main text. The projection-based methods seek the best cell orientation at each pixel in the design domain, but extra (and sometimes quite intensive) treatments to smoothly connect the neighbouring cells are needed. This is why the projection method \cite{Pantz_SIAM_JCO2008} is introduced. The reason, if viewed in a mathematical viewpoint, is a matter of total differential, i.e., knowing the partial derivatives at all pixel points does not mean one can naturally integrate them for a potential function, whose spatial gradients coincide with the prescribed partial derivatives. And the projection method is used to identify an approximate total differential form against the given partial derivatives. In contrast, the real and imaginary parts of a holomorphic function automatically satisfy the conditions for total differential, which are actually expressed by the Cauchy-Riemann equations. This is why no post-processing is needed for the proposed method based on conformal mapping.

The remainder of this article is organised as follows. In Section~\ref{sec: AHTO_plus_review}, the AABH-plus-based framework is briefly reviewed. This is followed by the introduction of using conformal mapping technique to derive the related formulation, and a novel optimisation framework towards graded orthotopic porous structures is established in Section~\ref{sec: AHTO_plus_confomal}. Then the issues concerning numerical implementation and corresponding sensitivity of the present approach are discussed in Section~\ref{sec: Numerical_and_sensitivity}. In Section~\ref{sec: Numerical_examples}, numerical examples are presented to illustrate the performance of the proposed method, and conclusions and perspectives compose of Section~\ref{sec: Discussion} at the end of the article.

\section{Review of AABH-plus-based optimisation framework}
\label{sec: AHTO_plus_review}
In this section, the AABH-plus-based optimisation framework is outlined first, and then key issues for its improvement are discussed.

\subsection{AABH-plus-based optimisation framework}
\label{subsec:AHTO_plus_Outline}
The AABH-plus-based framework \cite{Zhu_JMPS2019} mainly consists of three parts: the multiscale topology description function, AABH plus approach, and corresponding optimisation formulation.
\subsubsection{Multiscale description of graded porous structures}
\begin{figure}[!ht]
  \centering
  \includegraphics[width=.7\textwidth]{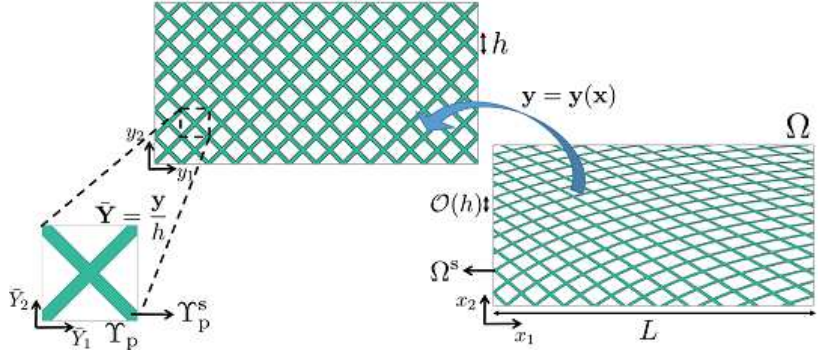}
  \caption{Multiscale representation of a graded porous structure generated through function composition given by Eq.~\eqref{Eq:Qusi_Peroidic_TDF}.}
  \label{fig:GMs_generation}
\end{figure}

The digital representation of a structure is usually carried out by a topology description function, say $\phi\left(\vx\right)$, defined by
\begin{equation}
  \label{Eq:TDF_of_GMs}
  \left\{
    \begin{aligned}
    &{\phi}({\vx}) > 0,&&\quad\forall \vx \in \Omega^\text{s},\\
    &{\phi}({\vx}) = 0,&&\quad\forall \vx \in \Gamma^\text{s},\\
    &{\phi}({\vx}) < 0,&&\quad\forall \vx \in \Omega\backslash \Omega^\text{s},
    \end{aligned}
  \right.
\end{equation}
where $\Omega^\text{s}$ and $\Omega$ denote the region occupied with solid materials and the whole design domain, respectively. A porous structure, as shown in Fig.~\ref{fig:GMs_generation}, is normally associated with two length scales: a microscale length $h$ on which microstructures are resolved and a macro length $L$ on which the design domain is measured, and $h \ll L$. Thus representing a porous structure pixel by pixel is costly. Here we can adopt a multiscale formulation to overcome this issue.

As shown in Fig.~\ref{fig:GMs_generation}, a smooth macroscopic mapping function $\vy = \vy(\vx)$ defined in the whole region $\Omega$ is introduced, so as to map the graded porous structure of interest to a periodic configuration of period $h$. Then the description of the resulting periodic structure can be confined within a single cell, say $\Upsilon_{\text{p}} = [0,1]^N$, which has been rescaled to be a unit cell, with $N$ being the dimensionality number. Suppose $\phi^{\text{p}}(\bar\vY)$ for $\bar{\vY} \in \Upsilon_\text{p}$ is adopted to identify the material layout in the unit cell. Then the actual graded porous configuration can be described through function composition, i.e.,
   \beq
   \label{Eq:Qusi_Peroidic_TDF}
   \phi(\vy(\vx)) = \phi^{\text{p}}\left(\frac{\mathbf{y}(\mathbf{x})}{h}\right) := \phi^{\text{p}}\left(\bar\vY(\vx)\right)
   \eeq
Here the structured unit cell $\Upsilon_\text{p}$ is termed as a ``matrix cell''. Eq.~\eqref{Eq:Qusi_Peroidic_TDF} effectively represents a porous structure by gradually rescaling, rotating and (iso-volumetrically) distorting the matrix cell in space. The onsite operation on the unit cell is captured by the locally Jacobian matrix defined by
   \beq
     \label{Eq:Jacobian}
    J_{ij} = \pd{y_i}{x_j}
   \eeq
     for $i,j = 1,\cdots, N.$

Eq.~\eqref{Eq:Qusi_Peroidic_TDF} makes use of the self-similar feature of graded microstructures. Hence fine mesh is only needed for resolving $\phi^\text{p}\left(\cdot\right)$, while the mapping function $\vy\left(\vx\right)$ is digitalised on the coarse-grained scale. Note that Eq.~\eqref{Eq:Qusi_Peroidic_TDF} can be generalised for cases with multiple matrix cells, and more details can be found in Xue et al. \cite{Xue_CMME2020}.
\subsubsection{AABH plus formulation}
\label{subsubsec:AABH_plus_approach}
The AABH plus formulation seeks for the scale-separated form of the multiscale formulation defined directly in the porous region, as the parameter
\beq
\label{Eq:Small_parameter}
\eps = \frac{h}{l}\rightarrow 0.
\eeq

For a graded porous region, the equilibrium equation defined in it should read
\beq
\label{Eq:Equilibrium_equation}
{-\pd{}{x_j} \left(\tilde{\mathbb{C}}^{\eps}_{ijkl}\left(\vy(\vx)\right) \pd{u^{\eps}_k}{x_l}\right) = f_i}
\eeq
for $i,j,k,l = 1,\cdots,N$, where $\vf = (f_i)$ is the body force per volume; $\vu^{\eps} = (u^{\eps}_i)$ denotes the displacement field; referring to Eq.~\eqref{Eq:Qusi_Peroidic_TDF}, $\tilde{\bm{\mathbb{C}}}^{\eps}=(\tilde{\mathbb{C}}^{\eps}_{ijkl}) $ is given by
\beq
\label{C_ijkl_extension}
\tilde{\mathbb{C}}^{\eps}_{ijkl}\left(\vy(\vx)\right) = \tilde{\mathbb{C}}^{\text{p}}_{ijkl}\left(\bar\vY(\vx)\right) =\left\{\begin{aligned}
  & \mathbb{C}_{ijkl}, \quad &&\vx\text{ in }\Omega^{\text{s}},\\
  & 0, &&\vx\text{ in }\Omega\backslash\Omega^{\text{s}},
  \end{aligned}\right.
\eeq
representing the elasticity tensor field; $\mathbb{C} = \left(\mathbb{C}_{ijkl}\right)$ is the elasticity tensor of the base material; a superscript $\eps$ is affiliated with a physical quantity to indicate it resolves the microstructural details in the design domain.

Given the structural periodicity when measured in $\bar \vY$ coordinates, as illustrated in Fig.~\ref{fig:GMs_generation}, the multiscale displacement field $\vu^\eps$ can be approximated by a series in terms of $\eps$, which is formed by a set of scale-separated function $\vu^{(i)}(\vx,\bar\vY)$ with $\vx$ now only measuring mean-field variations. Further asymptotic analysis leads to \cite{Zhu_JMPS2019}
\begin{gather}
\label{Eq:Homogenized_u}
{u}_i^{(0)}({\vx},\bar{\vY}) = {u}_i^{(0)}({\vx}) = {u}_i^{\text{H}},\\
\label{Eq:One_order_collector}
J_{mj}\pd{}{\bar{Y}_m} \left({\tilde{\mathbb{C}}^{\text{p}}}_{ijkl}  J_{nl}\pd{\xi_k^{st} }{\bar{Y}_n}\right) =J_{mj}\pd{{\tilde{\mathbb{C}}^{\text{p}}}_{ijkl} }{\bar{Y}_m},\\
\label{Eq:stiffness_equivalent}
\mathbb{C}_{ijkl}^{\text{H}} = \mathbb{C}_{ijkl}\cdot\left|\Upsilon_{\text{p}}^{\text{s}}\right| - \mathbb{C}_{ijst}J_{nt} \cdot\int_{\Upsilon_{\text{p}}^\text{s}} \pd{\xi_s^{kl}}{\bar{Y}_n} \dif\bar{\vY},
\end{gather}
where ${\vu}^{\text{H}} = \left(u_i^{\text{H}}\right)$ is the homogenised displacement field; the third-order tensor $\bm{\xi} = (\xi_k^{st})$ is referred to as a first-order corrector satisfying Eq.~\eqref{Eq:One_order_collector}, which is an artificially introduced intermediate quantity for computing the unit cell equivalent property; $\bm{\mathbb{C}^{\text{H}}} = (\mathbb{C}^{\text{H}}_{ijkl})$ denotes the homogenised elasticity tensor; $\left|\Upsilon_{\text{p}}^{\text{s}}\right|$ is given by
\beq
\label{Eq:Soild_volume_cell}
\left|\Upsilon_{\text{p}}^{\text{s}}\right| = \frac{\int_{\Upsilon_{\text{p}}^{\text{s}}}\dif\bar\vY} {\int_{\Upsilon_{\text{p}}}\dif\bar\vY},
\eeq
representing the volume fraction of solid materials in the unit cell. It is worth noting that since $\vJ$ depends on $\vx$, Eqs.~\eqref{Eq:One_order_collector} and \eqref{Eq:stiffness_equivalent} need to be solved for every $\vx \in \Omega$.

\subsubsection{Optimsation framework}
On taking the mapping function $\vy(\vx) $ as the macroscopic design variables, the TDF $\phi^{\text{p}}(\bar{\vY})$ of the matrix cell as the microscopic design variables, and the equivalent system compliance ${\mathcal{C}}^{\text{H}}$ as the target function, a compliance structures optimisation framework for graded porous configurations can be formulated as follows
\beq
\label{Eq:Optimsation_formulation}
\begin{aligned}
&\qquad\text{Find}\quad{\vy(\vx)} , \phi^{\text{p}}(\bar{\vY}),\\
&\qquad\text{Minimize} \quad\mathcal{C}^{\text{H}} =  \int_{\Omega} \mathbb{C}_{ijkl}^{\text{H}} \pd{u_i^{\text{H}}}{x_j} \pd{u_k^{\text{H}}}{x_l}\dif \vx,\\
&\text{s.t.}\\
&\qquad-\pd{}{{x}_j} \left({\mathbb{C}}_{ijkl}^{\text{H}} \pd{{u}_k^\text{H}}{{x}_l}\right) = f_i, \quad   \left.u_i^{\text{H}}\right|_{\Gamma_u} = 0,\quad \left.\mathbb{C}_{ijkl}^{\text{H}}n_j\pd{u_k^{\text{H}}}{x_l}\right|_{\Gamma_t} = t_i,\\
&\qquad J_{mj}\pd{}{\bar{Y}_m} \left(\tilde{\mathbb{C}}^{\text{p}}_{ijkl} J_{nl}\pd{\xi_k^{st}}{\bar{Y}_n}\right)  = J_{mj} \pd{\tilde{\mathbb{C}}^{\text{p}}_{ijst} }{\bar{Y}_m},\\
&\qquad \mathbb{C}_{ijkl}^{\text{H}} = \mathbb{C}_{ijkl}\cdot\left|\Upsilon_{\text{p}}^{\text{s}}\right| - \mathbb{C}_{ijst}J_{nt} \cdot\int_{\Upsilon_{\text{p}}^\text{s}} \pd{\xi_s^{kl}}{\bar{Y}_n} \dif\bar{\vY},\\
&\qquad V_f \le \bar{V},\\
\end{aligned}
\eeq
where $\Gamma_u$ and $\Gamma_t$ represent the boundary sections imposed with Dirichlet and Neumann conditions, respectively; $\vt = (t_i)$ denotes the surface force per area acting on the $\Gamma_t$; $\vn = (n_i)$ is the outer normal vector of the design domain; $V_f$ is the volume fraction of the region part with solid materials; $\bar{V}$ denotes the admissible upper bound of the volume fraction.

\subsection{Challenging issues}
Up to this stage, the potential of the AABH-plus-based framework has not been fully realised, mainly caused by the following reasons. Firstly, as pointed out earlier, since the Jacobian matrix depends on the macroscopic coordinate $\vx$, the first-order corrector $\boldsymbol{\xi}$ should be solved for point-wisely. This is computationally expensive. To address this issue, actions were taken \cite{Zhu_JMPS2019, Xue_Smo2020}, but improvement is still in demand. Secondly, the mapping function $\vy = \vy\left(\vx\right)$ is up to now restricted to polynomials, and it appears that the optimisation result is quite sensitive to the choice of basis functions. For compliance optimisation, high distortion of the matrix cell is not favoured, and the corresponding treatments for avoiding that are desired.

The present article is aimed to demonstrate that the use of conformal mapping appropriately addresses the challenging issues mentioned above. With the excellent properties of conformal mapping, the proposed framework manages to archive a delicate balance among computational efficiency, gradual varying freedom, and limiting distortion, for compliance optimisation of two-dimensional porous configurations.

\section{AABH-plus-based optimisation framework combined with conformal mapping}

\label{sec: AHTO_plus_confomal}
In this section, we demonstrate how the concept of conformal mapping is integrated in the AABH-plus-based optimisation framework. Besides, the case of rectangular design domain where analytical solutions are available is discussed.

\subsection{Properties of graded porous structures generated through conformal mapping}
\label{subsec:Properties_GMs_conformal_mapping}

\subsubsection{Problem set-up}
\label{subsubsec:Problem_set_up}
\begin{figure}[!ht]
  \centering
  \includegraphics[width=.45\textwidth]{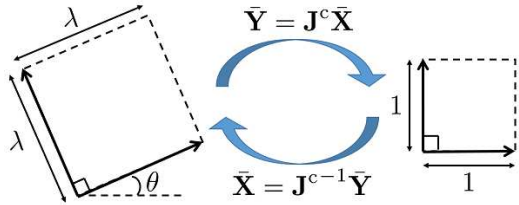}
  \caption{Mapping between cells.}
  \label{fig:Cells_Map}
\end{figure}

In the original AABH-plus-based framework \cite{Zhu_JMPS2019}, the matrix cell can be rescaled, rotated and distorted in accordance with the onsite macroscopic mapping function $\vy\left(\vx\right)$. Now we restrict the porous configuration of our attention to be generated through conformal mapping, where the angles about any two crossing microstructural members stay unchanged after mapping. Thus the deformation paradigms of microstructures are limited to rescaling and rotation only. Here we use the symbol $\vy^{\text{c}}\left(\vx\right)$ to emphasise that the corresponding mapping is conformal. Hence, the Jacobian matrix $\vJ^{\text{c}}$ of $\vy^{\text{c}}(\vx)$ can be specified to be
\beq
\label{Eq:J_R}
{\textbf{J}^{\text{c}}}=  \frac{1}{\lambda}\textbf{R}^\text{-1},
\eeq
where $\textbf{R}$ denotes a rotation matrix satisfying $\textbf{R}^\top\textbf{R} = \textbf{I}$ with $\textbf{I}$ being the identity matrix. For two-dimensional cases, $\vJ^{\text{c}}$ can be further identified by
\beq
\label{Eq:J_form}
{\textbf{J}^{\text{c}}}=  \frac{1}{\lambda}\begin{pmatrix} \cos \theta &\sin\theta \\ -\sin\theta & \cos\theta \end{pmatrix},
\eeq
where $\lambda$ and $\theta$ denote the scaling factor and the rotation angle of the matrix cell in a counterclockwise sense, respectively, as shown in Fig.~\ref{fig:Cells_Map}. Here $\theta = 0$ stands for the situation when one edge of the matrix cell is parallel to the $\bar Y_1$-axis. Note that both $\lambda$ and $\theta$ are functions of the macroscopic coordinate $\vx$, where $\lambda>0$ holds for $\vx \in \Omega$.

Based on Eq.~\eqref{Eq:J_form}, the components of the resulting Jacobian matrix $\vJ^\text{c}$, as defined by Eq.~\eqref{Eq:Jacobian}, should be linked with the mapping function $\vy^\text{c}\left(\vx\right)$ by
\begin{subequations}
  \label{Eq:Diff_y}
  \begin{alignat}{2}
    \label{Eq:Diff_y1}
    &\dif y^\text{c}_1 =\frac{1}{\lambda}\cos \theta\mathrm{d}x_1 +\frac{1}{\lambda} \sin\theta \mathrm{d}x_2;\\
    \label{Eq:Diff_y2}
    &\dif y^\text{c}_2 =-\frac{1}{\lambda}\sin \theta\mathrm{d}x_1 +\frac{1}{\lambda} \cos\theta \mathrm{d} x_2.
\end{alignat}
\end{subequations}
To ensure the existence of the mapping function $\vy^{\text{c}}(\vx)$, $\dif y_1^{\text{c}}$ and $\dif y_2^{\text{c}}$ must take exact differential forms, indicating that the components of $\vJ^{\text{c}}$ should be interrelated by the conditions for total differential, i.e.,
\begin{subequations}\label{Eq:Integrable}
\begin{alignat}{2}
  \label{Eq:y1_Integrable}
  &\pd{\left(\frac{1}{\lambda}\cos\theta\right)}{x_2}= \pd{\left(\frac{1}{\lambda}\sin\theta\right)}{x_1};\\
  \label{Eq:y2_Integrable}
  &\pd{\left(-\frac{1}{\lambda}\sin\theta\right)}{x_2}= \pd{\left(\frac{1}{\lambda}\cos\theta\right)}{x_1}.
\end{alignat}
\end{subequations}
It is shown in Appendix~\ref{Appendix:Cauchy_Riemann} that Eqs.~\eqref{Eq:Integrable} lead to
\begin{subequations}\label{Eq:Cauchy_Riemann}
\begin{alignat}{2}
    \label{Eq:Cauchy_Riemann_1}
    &\pd{\left(\ln\lambda\right)}{x_1}=\pd{\theta}{x_2};\\
    \label{Eq:Cauchy_Riemann_2}
    &\pd{\left(\ln\lambda\right)}{x_2}=-\pd{\theta}{x_1},
\end{alignat}
\end{subequations}
which means that $\ln\lambda$ and $\theta$ satisfy the Cauchy-Riemann equations.

To this end, if we introduce a complex argument
\beq
\label{Eq: Complex argument}
z = x_1 + i x_2,
\eeq
where $i$ represents the imaginary unit, $\ln\lambda$ is conjugate to $\theta$, i.e., they are just the real and the imaginary parts of a holomorphic function. As one conclusion from a holomorphic function, both $\ln\lambda$ and $\theta$ defined in $\Omega$ are harmonic, that is,
\begin{subequations}\label{Eq:Harmonic_equation}
  \begin{alignat}{2}
      \label{Eq:Harmonic_lnlmd}
      &\nabla^2 \left(\ln\lambda\right) = 0;\\
      \label{Eq:Harmonic_theta}
      &\nabla^2 \theta = 0.
  \end{alignat}
  \end{subequations}
\subsubsection{Boundary conditions}
\label{subsubsec:Boundary_condition_Cauchy_Riemann}
\begin{figure}[!ht]
  \centering
  \includegraphics[width=.3\textwidth]{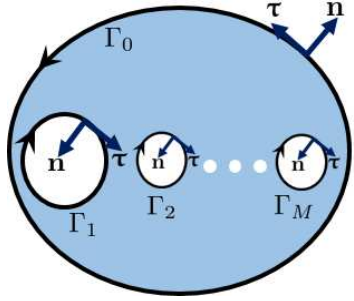}
  \caption{Domain of definition. The design domain may be multiply connected.}
  \label{fig:Cauchy_Riemann_Difinite_Domain}
\end{figure}

Now we consider imposing proper boundary conditions for the two harmonic functions, $\ln\lambda$ and $\theta$. The boundary of a domain, as shown in Fig.~\ref{fig:Cauchy_Riemann_Difinite_Domain}, may consist of several simple closed oriented curves, denoted by $\Gamma = \cup_{i= 0}^{M}\Gamma_{i}$, where $\Gamma_0$ is set to be the outer boundary of the region; $\Gamma_i, i = 1,\cdots,M$ form the boundaries of its $M$ inner "holes"; $\vn$ and $\bm\uptau$ represent the outer normal vector and the tangent vector of the boundary, respectively. Here the direction of $\bm\uptau$ should be specified for a multiply connected region. Here $\bm\uptau$ is set counterclockwise on the outer boundary $\Gamma_0$, and clockwise on all other parts, i.e. $\Gamma_i, i=1,\cdots,M$, as shown in Fig.~\ref{fig:Cauchy_Riemann_Difinite_Domain}.

Now we demonstrate that one only needs to impose boundary conditions for either $\ln\lambda$ or $\theta$. Hence both $\ln\lambda$ and $\theta$ are almost fully determined (up to a constant) in the multiply domain definition. Here we assume the values of $\ln\lambda$ are given on $\Gamma$. The directional derivative of $\ln\lambda$ along the tangential direction on the region boundary, which is also known, is expressed by
\beq
\label{Eq:Directional_derivative}
\pd{\left(\ln\lambda\right)}{\bm\uptau}= \pd{(\ln\lambda)}{x_1}\cos(\bm\uptau,x_1) + \pd{(\ln\lambda)}{x_2}\cos(\bm\uptau,x_2),
\eeq
where
\begin{subequations}\label{Eq:Tangential_normal}
  \begin{alignat}{2}
      \label{Eq:Tangential_normal_1}
      &\cos(\vn,x_1) = \cos(\bm\uptau,x_2);\\
      \label{Eq:Tangential_normal_2}
      &\cos(\vn,x_2) = -\cos(\bm\uptau,x_1).
  \end{alignat}
  \end{subequations}
Incorporating Eqs.~\eqref{Eq:Cauchy_Riemann} and ~\eqref{Eq:Tangential_normal} into Eq.~\eqref{Eq:Directional_derivative} gives
\beq
\label{Eq:Normal_derivative}
\pd{\theta}{\vn} = -\pd{\left(\ln\lambda\right)}{\bm\uptau},
\eeq
which indicates that the boundary information for $\theta$ is  also known by means of Neumann boundary condition.

Therefore, a pair of definite problems for $\ln\lambda$ and $\theta$ are set to be
\begin{subequations}
\label{Eq:BC_determine_solution}
\begin{alignat}{2}
  \label{Eq:BC_determine_solution_lnlmd}
  &\nabla^2\left(\ln\lambda\right)=0,&\quad &\left(\ln\lambda\right)|_{\Gamma}=\ln{\lambda_b};    \\
  \label{Eq:BC_determine_solution_theta}
  &\nabla^2\theta=0,&\quad &\pd{\theta}{\vn}\bigg|_{\Gamma}=-\pd{\ln{\lambda_b}}{\bm \uptau}.
\end{alignat}
\end{subequations}
where $\lambda_b$ denotes the values of $\lambda$ on the boundary. Note that, if $\theta$ is one solution to problem~\eqref{Eq:BC_determine_solution_theta}, so is $\theta+c$, with $c$ being an arbitrary constant. Hence adding a supplementary condition
\beq
\label{Eq:Supplementary_constant}
\int_{\Omega}\theta\dif\vx =\bar\theta\int_{\Omega}\dif\vx=c{\int_{\Omega}\dif\vx}
\eeq
guarantees the uniqueness of the solution to problem~\eqref{Eq:BC_determine_solution_theta}, where $\bar\theta$ is a constant representing the mean value of the rotation angle.

Up to this stage, we manage to show that $\ln\lambda$ and $\theta$ within the design domain are fully determined with the boundary values $\ln\lambda_b$ as well as the mean notation angle $\bar\theta$. Hence we are enabled to just confine the design variables as $\ln\lambda_b$ on a collection of boundary points, as well as the averaged rotational angle $\bar\theta$.

For a given design domain $\Omega$, one may further solve Eqs.~\eqref{Eq:BC_determine_solution}~--~\eqref{Eq:Supplementary_constant} based on the method of Green's functions and obtain
\begin{subequations}
\label{Eq:Cauchy_Riemann_solution}
\begin{alignat}{2}
\label{Eq:lnlmd_solution}
&{\ln\lambda\left(\vx\right)=\int_{\Gamma}\ln{\lambda}\left(\vx'\right)\cdot\pd{{G}^{\textup{D}}\left(\vx,\vx'\right)}{\vn}  \dif S}, \quad \vx' \in \Gamma;\\
\label{Eq:Theta_solution}
&\theta\left(\vx\right)=-\int_{\Gamma}G^{\textup{N}}\left(\vx,\vx'\right)\cdot\pd{\theta\left(\vx'\right)}{\vn}\dif S + \bar\theta,\quad \vx' \in \Gamma,
\end{alignat}
\end{subequations}
where ${G}^{\textup{D}}\left(\vx,\vx'\right)$ and $G^{\textup{N}}\left(\vx,\vx'\right)$ are the Green's functions corresponding to Dirichlet and Neumann conditions imposed for Laplace's equation, satisfying
\begin{subequations}
  \label{Eq:Green_function}
    \begin{alignat}{2}
      \label{Eq:Green_function_dirichlet}
        &\nabla^2 G^{\textup{D}}\left(\vx,\vx'\right) = \delta\left(\vx-\vx'\right),& \quad & G^{\textup{D}}\left(\vx,\vx'\right)\vert_{\Gamma}  = 0; \\
      \label{Eq:Green_function_neumann}
        &\nabla^2 G^{\textup{N}}\left(\vx,\vx'\right) = \delta\left(\vx-\vx'\right) -\frac1{\int_{\Omega}\dif \vx'},& \quad &\pd{G^{\textup{N}}\left(\vx,\vx'\right)}{\vn} \bigg|_{\Gamma} = 0,
    \end{alignat}
\end{subequations}
where $\delta\left(\cdot\right)$ is the Dirac-delta function.

It is worth noting that both ${G}^{\textup{D}}\left(\vx,\vx'\right)$ and $G^{\textup{N}}\left(\vx,\vx'\right)$ are independent of boundary conditions and just rely on the shape of $\Omega$. Hence we can get Green's functions in an offline stage. When the boundary values are changed, one simply inserts them into Eq.~\eqref{Eq:Cauchy_Riemann_solution}, and $\ln\lambda$ and $\theta$ can be obtained without solving Laplace's equation repeatedly.

\subsubsection{Determination of the mapping function}
\label{Determine_mapping}
To fully represent a porous configuration with Eq.~\eqref{Eq:Qusi_Peroidic_TDF}, one still needs to express the mapping function $\vy^{\text{c}}(\vx)$. If $\Omega$ is a simply connected region, it can be seen from Eq.~\eqref{Eq:Integrable} that the line integrals of $\dif y^\text{c}_1$ and $\dif y^\text{c}_2$ defined by Eq.~\eqref{Eq:Diff_y} in $\Omega$ are path independent. Therefore $\vy^{\text{c}}(\vx)$ can be expressed by
\begin{subequations}
\label{Eq:Integral_y}
\begin{alignat}{2}
  \label{Eq:Integral_y1}
    &y_1^\text{c}(\vx) = \int_{\vx_0}^{\vx}\frac{1}{\lambda}\cos \theta\mathrm{d}x_1 +\frac{1}{\lambda} \sin\theta \mathrm{d}x_2 + c_1, \\
  \label{Eq:Integral_y2}
    &y_2^\text{c}(\vx) = \int_{\vx_0}^{\vx}-\frac{1}{\lambda}\sin \theta\mathrm{d}x_1 +\frac{1}{\lambda} \cos\theta \mathrm{d} x_2 +c_2,
\end{alignat}
\end{subequations}
where $c_1$ and $c_2$ are the values of $y^\text{c}_1$ and $y^{\text{c}}_2$ at $\vx_0$, respectively, which cast no effects on the configuration of the resulting graded porous structure because of the scale separation.

If $\Omega$ is a multiply connected region as shown in Fig.~\ref{fig:Cauchy_Riemann_Difinite_Domain}, in addition to Eq.~\eqref{Eq:Integrable}, the following conditions
\begin{subequations}
  \label{Eq:Multi_connected}
  \begin{alignat}{2}
    \label{Eq:Multi_connected_y1}
      &\oint_{\Gamma_{i}}\frac{1}{\lambda}\cos \theta\mathrm{d}x_1 + \frac{1}{\lambda} \sin\theta \mathrm{d}x_2 = 0,\quad{i = 1,\cdots ,M}; \\
    \label{Eq:Multi_connected_y2}
      & \oint_{\Gamma_{i}}-\frac{1}{\lambda}\sin \theta\mathrm{d}x_1 +\frac{1}{\lambda} \cos\theta \mathrm{d} x_2=0,\quad{i = 1,\cdots ,M,}
  \end{alignat}
  \end{subequations}
should also hold so to ensure the integrability of $\dif y_1^{\text{c}}$ and $\dif y_2^{\text{c}}$.
\subsection{Microstructural size control}
\label{Extremal_principle}
As $\ln\lambda$ is harmonic in $\Omega$, the maximum principle reads
\begin{subequations}\label{Eq:Extremal_principle}
  \begin{alignat}{2}
    \label{Eq:Extremal_principle_max}
    &\max_{\vx\in\Omega}\ln\lambda=\max_{\vx \in\Gamma} \ln{\lambda_b} ;\\
    \label{Eq:Extremal_principle_min}
    &\min_{\vx\in\Omega}\ln\lambda=\min_{\vx \in\Gamma} \ln{\lambda_b} .
  \end{alignat}
\end{subequations}
Therefore, the range for $\lambda$ in the domain can be fully determined by simply controlling the values of $\lambda$ on the boundaries, i.e. $\lambda_b$. Consequently, the upper and lower bounds of the size of microstructural members constituting a porous configuration can be explicitly imposed through the design variable, $\ln\lambda$. This is of practical value when porous configurations are processed with additive manufacturing techniques, which often yield a minimal printable size. To actually impose the constraint, we first use $d_{\min}$ to denote the (non-dimensioned) minimal size of the microstructures member in the matrix cell, with reference to Eq.~\eqref{Eq:Extremal_principle_min}, the minimal size $d_{\min}$ of the entire structure can be expressed by
\beq
\label{Eq:Minimal_size}
d_{\min} = h D_{\min}\cdot {\min_{\vx \in\Gamma} {\lambda_b}},
\eeq
where $h$ is recalled to be the characteristic length of the unit cells as indicated in Fig.~\ref{fig:GMs_generation}. Therefore, as long as the minimum of the $\lambda_b$ on $\Gamma$ satisfies
\beq
\label{Eq:Size_constraint}
\min_{\vx \in\Gamma}{\lambda_b} \ge \frac{p_{\min}}{hD_{\min}},
\eeq
it can be ensured that the minimal size of the microstructure components falls be low the minimal printable size denoted by $p_{\min}$.

\subsection{Homogenisation}
Compared with general cases, the porous structures considered here result only from the rescaling and rotation of the matrix cell. This also helps reduce the computational cost for homogenisation, if compared with the general case briefly reviewed in Sec.~\ref{subsubsec:AABH_plus_approach}.

To see this, we start with the computation of the elasticity tensor homogenised from the (unit) matrix cell. To this end, we let $\textbf{J} = \textbf{I}$ in Eqs.~\eqref{Eq:One_order_collector} and~\eqref{Eq:stiffness_equivalent} with $\textbf{I}$ being the identity matrix, and obtain
\begin{gather}
  \label{Eq:Represent_cell_one_order_collector}
  \pd{}{{\bar{Y}}_j} \left(\tilde{\mathbb{C}}^{\text{p}}_{ijkl} \pd{{\hat\xi_k^{st}}}{{\bar{Y}}_l}\right)= \pd{\tilde{\mathbb{C}}^{\text{p}}_{ijst} }{\bar{Y}_j};\\
  \label{Eq:Represent_stiffness_equivalent}
  {\hat{\mathbb{C}}_{ijkl}^{\text{H}}}
= {{\mathbb{C}}}_{ijkl}\cdot |\Upsilon_{\text{p}}^\text{s}| - {{\mathbb{C}}}_{ijst}\cdot\int_{\Upsilon_{\text{p}}^{\text{s}}}   \pd{   {\hat\xi_s^{kl}}   }{\bar{Y}_t} \dif \bar\vY,
\end{gather}
where $\hat{\bm{\xi}}$ and $\hat{\bm{\mathbb{C}}}$ represent the first-order corrector and homogenised elasticity tensor of the matrix cell when $\theta = 0$ and $\lambda=1$, respectively.

It is proved in Appendix~\ref{Appendix:One_order_corrector} that when the base materials are elastically isotropic, the resulting first-order corrector $\bm{\xi}$ and the homogenised elasticity tensor $\mathbb{C}^{\text{H}}$ for arbitrary $\lambda$ and $\theta$ can be linked with that of the unit matrix cell by
\beq
\label{Eq:Conformal_cell_one_order_collector}
{\xi_w^{uv}}={\lambda}R_{us}R_{vt}R_{wk}\hat\xi_k^{st}
\eeq
and
\beq
\label{Eq:Conformal_cell_stiffness}
{\mathbb{C}_{ijkl}^{\text{H}}}=R_{ip}R_{jq}R_{ks}R_{lt}{\hat{\mathbb{C}}_{pqst}^{\text{H}}},
\eeq
respectively.

Eq.~\eqref{Eq:Conformal_cell_stiffness} is intuitive. Due to the conformality, each constituting cell of the porous structure is effectively obtained by rescaling and rotating the matrix cell. Furthermore, the equivalent properties of unit cells are independent of the scaling factor due to scale separation. Hence the equivalent elasticity tensor $\mathbb{C}^{\text{H}}$ of each unit cell is determined by rotating the initial homogenised elasticity tensor $\hat{\mathbb{C}}^{\text{H}}$ counterclockwise with angle $\theta$.

In the present case, therefore, one simply needs to calculate the homogenised elasticity tensor for a benchmarked case with $\lambda = 1$ and $\theta = 0$. This significantly reduces the computational cost for general cases where the cell problem has to be resolved point by point.

\subsection{Optimisation}
Based on the discussion above, an AABH-plus-based optimisation framework combined with conformal mapping is formulated to be
\beq
\label{Eq:Optimation_formula_Conformal_GMs}
\begin{aligned}
&\qquad\text{Find}\quad{\ln{{\lambda_b}},\bar\theta,\phi^{\text{p}}(\bar{\vY})},\\
&\qquad\text{Minimize} \quad\mathcal{C}^{\text{H}} =   \int_{\Omega} \mathbb{C}_{ijkl}^{\text{H}} \pd{u_i^{\text{H}}}{x_j} \pd{u_k^{\text{H}}}{x_l}\dif \vx,\\
&\text{s.t.}\\
&\qquad-\pd{}{{x}_j} \left({\mathbb{C}}_{ijkl}^{\text{H}} \pd{{u}_k^\text{H}}{{x}_l}\right) = f_i, \quad\left.u_i^{\text{H}}\right|_{\Gamma_1^{\text{H}}} = 0,\quad \left.\mathbb{C}_{ijkl}^{\text{H}}n_j\pd{u_k^{\text{H}}}{x_l}\right|_{\Gamma_2^{\text{H}}} = t_i,\\
&\qquad\nabla^2\left(\ln\lambda\right)=0,\quad \left(\ln\lambda\right)|_{\Gamma}=\ln{\lambda_b},\\
&\qquad\nabla^2\theta=0,\quad \pd{\theta}{\vn}\bigg|_{\Gamma}=-\pd{\ln{\lambda_b}}{\bm \uptau},\quad \int_{\Omega}\theta\dif\vx =  \bar\theta{\int_{\Omega}\dif\vx},\\
&\qquad\pd{}{\bar{Y}_j} \left(\tilde{\mathbb{C}}^{\text{p}}_{ijkl} \pd{{\hat{\xi}_k^{st}}}{\bar{{Y}}_l}\right)= \pd{\tilde{\mathbb{C}}^{\text{p}}_{ijst} }{\bar{Y}_j},\\
&\qquad{\hat{\mathbb{C}}_{ijkl}^{\text{H}}}
= {{\mathbb{C}}}_{ijkl}\cdot |\Upsilon_{\text{p}}^\text{s}|- \int_{\Upsilon_{\text{p}}^{\text{s}}} {{\mathbb{C}}}_{ijst}  \pd{   {\hat{\xi}_s^{kl}}   }{\bar{Y}_t} \dif \bar{\vY},\\
&\qquad{\mathbb{C}_{ijkl}^{\text{H}}}=R_{ip}R_{jq}R_{ks}R_{lt}{\hat{\mathbb{C}}_{pqst}^{\text{H}}}, \\
&\qquad \min_{\vx \in\Gamma}{\lambda_b} \ge \frac{p_{\min}}{hD_{\min}},\\
&\qquad V_f \le \bar{V},\\
\end{aligned}
\eeq
where $V_f$ represents the volume fraction of the solid parts constituting the porous configuration.

Note that we choose $\ln\lambda_b$, rather than the rotational angle $\theta$, as the microscopic design variables. This is because $\theta$ appears in the optimisation formulation by means of $\sin\theta$ or $\cos\theta$. Whichever the case, ambiguities caused by a jump of $\pi$ in $\theta$ may appear. Note that this is a problem reported in other frameworks \cite{Groen_IJNME2018}. Such ambiguities can be circumvented by employing $\ln\lambda_b$ as the design variables.

\subsection{Special cases of rectangular design domain}
\begin{figure}[!ht]
  \centering
  \includegraphics[width=.6\textwidth]{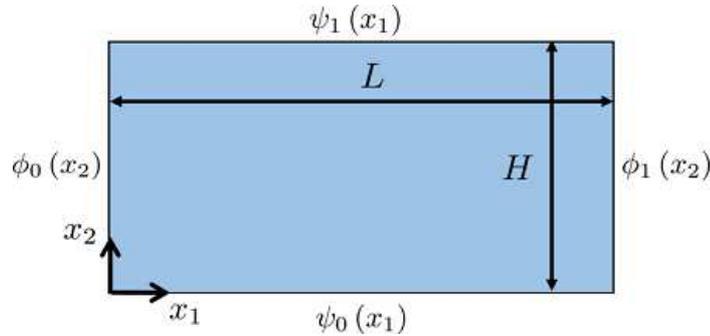}
  \caption{The case where the design domain takes a rectangular shape.}\label{fig:Rectangle_domian}
\end{figure}

When the design domain takes a rectangular shape, where Laplace's equation may be analytically solved, the computational burden for conducting optimisation here can be further levied. To see this, we consider a design domain as shown in Fig.~\ref{fig:Rectangle_domian}, where $\Omega = [0,L]\times[0,H]$; $\psi_0\left(x_1\right),\psi_1\left(x_1\right),\phi_0\left(x_2\right)$ and $\phi_1\left(x_2\right)$ represent the values of $\ln\lambda$ correspondingly on $\Gamma$. Note that, the compatibility conditions, $\phi_0\left(0\right)=\psi_0\left(0\right),\phi_0\left(H\right)=\psi_1\left(0\right),\phi_1\left(0\right)=\psi_0\left(L\right)$, and $\phi_1\left(H\right)=\psi_1\left(L\right)$ should also hold.

Hence the governing equations of $\ln\lambda$ and $\theta$ outlined in the optimisation formulation~\eqref{Eq:Optimation_formula_Conformal_GMs}
can be re-written by
\begin{subequations}
\label{Eq:BC_determine_solution_rectangle}
\beq
\label{Eq:Rectangle_lnlmd}
\left\{
\begin{aligned}
&\nabla^2\left(\ln\lambda\right)=0,    \quad  \left(0 \le x_1 \le L, \text{ } 0 \le x_2 \le H\right),\\
&\left(\ln\lambda\right)|_{x_1=0}=\phi_0\left(x_2\right),\quad\left(\ln\lambda\right)|_{x_1=L}=\phi_1\left(x_2\right), \\
&\left(\ln\lambda\right)|_{x_2=0}=\psi_0\left(x_1\right),\quad\left(\ln\lambda\right)|_{x_2=H}=\psi_1\left(x_1\right), \\
\end{aligned}
\right.
\eeq\\
\text{and}
\beq
\label{Eq:Rectangle_theta}
\:\;\;
\left\{
\begin{aligned}
&\nabla^2\theta=0,\quad  \left(0 \le x_1 \le L, \text{ } 0 \le x_2 \le H\right),\\
&\left.\pd{\theta}{x_1}\right|_{x_1=0}=-\pd{\phi_0\left(x_2\right)}{x_2},\quad\left.\pd{\theta}{x_1}\right|_{x_1=L} =-\pd{\phi_1\left(x_2\right)}{x_2}, \\
&\left.\pd{\theta}{x_2}\right|_{x_2=0}=\pd{\psi_0\left(x_1\right)}{x_1}, \:\:\;\quad\left.\pd{\theta}{x_2}\right|_{x_2=H}=\pd{\psi_1\left(x_1\right)}{x_1},\\
&{\int_{\Omega}\theta\dif\vx} = \bar\theta{\int_{\Omega}\dif\vx},
\end{aligned}
\right.
\eeq
\end{subequations}
respectively.

With the use of the method of separation of variables, the solutions to problems~\eqref{Eq:BC_determine_solution_rectangle} can be expressed in terms of series given by
\begin{subequations}
  \label{Eq:Rectangle_solution}
  \beq
  \begin{footnotesize}
  \label{Eq:Rectangle_lnlmd_solution}
  \begin{aligned}
  \ln\lambda =&\sum_{k=1}^\infty \left(a_k\cosh\frac{k\pi{x_2}}{L}+b_k\sinh\frac{k\pi{x_2}}{L}\right)\sin\frac{k\pi{x_1}}{L}+\left(c_k\cosh \frac{k\pi{x_1}}{H}+d_k\sinh\frac{k\pi{x_1}}{H}\right)\sin\frac{k\pi{x_2}}{H}\\
    &+\left(\psi_0\left(0\right)-\psi_0\left(L\right)+\psi_1\left(L\right)-\psi_1\left(0\right)\right)\frac{x_1x_2}{HL}
  +\left( \psi_0\left(L\right)-\psi_0\left(0\right) \right)\frac{x_1}{L}\\
  &+\left( \psi_1\left(0\right)-\psi_0\left(0\right) \right)\frac{x_2}{H}+\psi_0\left(0\right),
  \end{aligned}
  \end{footnotesize}
  \eeq
  and
  \beq
  \begin{footnotesize}
  \label{Eq:Rectangle_theta_solution}
  \begin{aligned}
    \theta =&\sum_{k=1}^\infty \left(a_k\sinh\frac{k\pi{x_2}}{L}+b_k\cosh\frac{k\pi{x_2}}{L}\right)\cos\frac{k\pi{x_1}}{L}- \left(c_k\sinh\frac{k\pi{x_1}}{H}+d_k\cosh\frac{k\pi{x_1}}{H}\right)\cos\frac{k\pi{x_2}}{H}\\
    &+\left(\psi_0\left(0\right)-\psi_0\left(L\right)+\psi_1\left(L\right)-\psi_1\left(0\right)\right)\frac{3 \left(x_2\right)^2-3\left(x_1\right)^2+H^2-L^2}{6HL}\\
    &+\left( \psi_0\left(L\right)-\psi_0\left(0\right) \right)\frac{x_2-\frac{H}{2}}{L} -\left( \psi_1\left(0\right)-\psi_0\left(0\right) \right)\frac{x_1-\frac{L}{2}}{H} +c,
  \end{aligned}
  \end{footnotesize}
  \eeq
\end{subequations}
where the detailed expressions for $a_k$,$b_k$,$c_k$, and $d_k$, for $k = 1,2,\cdots,$ as well as the derivations for Eqs.~\eqref{Eq:Rectangle_solution} can be found in Appendix~\ref{Appendix:Rectangular_solution_Cauchy_Riemann}. Therefore the procedure of resolving Laplace's equation can be skipped at each step during optimisation. Instead, one may just insert the design variables which are now $\psi_0\left(x_1\right),\psi_1\left(x_1\right),\phi_0\left(x_2\right)$ and $\phi_1\left(x_2\right)$, into the expressions~\eqref{Eq:Rectangle_solution}.
\section{Issues on numerical implementation}
\label{sec: Numerical_and_sensitivity}
In this section, several key numerical issues on implementing the present optimisation method, such as the digital representation of design variables, sensitivity analysis, are considered.

\subsection{Digital representation of design variables}
\label{subsec:Parametric_design_variables}
The design variables of the proposed optimisation framework are the logarithm of the scaling factor evaluated on the domain boundary, that is, $\ln\lambda_b$ and the TDF $\phi\left(\bar\vY\right)$ of the matrix cell. Now we need to determine appropriate digital representation for them.

\subsubsection{Microscopic design variables}
\label{subsubsec:Control_parameter_cell}
Following the preceding studies \cite{Zhu_JMPS2019, Xue_Smo2020}, the MMC framework \cite{Guo_JAM2014, Zhang_SMO2015, Zhang_CM2016}, is employed to represent the material layout within the matrix cell. In MMC framework, $\phi^\text{p}(\bar\vY)$ can be represented by the geometric parameters of components, and the microscopic design variables can be collected as
\beq
\label{Eq:Micro_design_variables}
\mathbf{D} = {\left(\bar{Y}^1_{01},\bar{Y}^1_{02}, a^1,b^1, \alpha^1, \cdots, \bar{Y}^n_{01},\bar{Y}^n_{02}, a^n,b^n, \alpha^n\right)}^{\top},
\eeq
where $(\bar{Y}^i_{01}, \bar{Y}^i_{02})$, $a^i$, $b^i$ and $\alpha^i$ denote the central coordinates, the half-length, the half-width, and the oriented angle of $i$th moving morphable component, respectively.
\subsubsection{Macroscopic design variables}
\label{subsubsec:Control_parameter_mapping}
In this paper, interpolation is used to describe $\ln\lambda_b$. Here we assume that the boundary of the design domain always consists of simple closed curves, each of which can be parameterised by the arc-length variable. Mathematically, this can be given by
 \beq
 \label{Eq:Curve_parameter}
 \Gamma_i = \Gamma_i\left(s_i\right), \quad s_i \in \left[0,l_i\right]
 \eeq
for $i = 0,\cdots,M$, where $s_i$ and $l_i$ are the arc length and the perimeter of $\Gamma_i$, respectively. Furthermore, $\ln\lambda_b$ defined on the boundary can be expressed by
\beq
\label{Eq:Curve_lnlmd}
\ln\lambda_b = \ln\lambda_b(\Gamma_i\left(s_i\right)),
\eeq
where $\ln\lambda_b(\Gamma_i(0)) = \ln\lambda_b(\Gamma_i(l_i))$ for continuity. Now the arc length parameter is then linked with the boundary condition along the boundary tangent, i.e.,
\beq
\label{Eq:Arc_length_parameterisation_directional_derivative}
\pd{\ln\lambda_b}{\bm{\uptau}} = \frac{\dif \ln\lambda_b}{\dif s_i},\quad \vx \in \Gamma_i.
\eeq

Select $(n_i+1)$ points $s_0^i$, $s_1^i$, $\cdots$, $s_{n_i}^i$ equidistantly on $\left[0,l_i\right]$, where $s_0^i = 0$ and $s_{n}^i = l_i$, and define $f^i_j$ as the values of $\ln\lambda_b\left(s_i\right)$ at point $s_j^i$, that is,
\beq
\label{Eq:Curve_interpolation_node}
f_j^i:=\ln\lambda_b\left(\Gamma_i\left(s^i_j\right)\right)
\eeq
for $j = 0,\cdots,n_i$, where $h^i = \frac{l_i}{n_i}$ denotes the arc-length difference between the neighbouring nodes on $\Gamma_i$.

Through interpolation, $\ln\lambda_b(\Gamma_i(s_i))$ is given by
\beq
\label{Eq:Curve_lnlmd_parameterisation}
\ln\lambda_b(\Gamma_i(s_i)) =\sum^{n_i}_{j=0}\omega_j\left(s_i\right)f^i_j ,\quad   s_i \in \left[0,l_i\right]
\eeq
for $j = 0, \cdots, n_i$, where the $f^i_j$ are interpolation nodes on $\Gamma_i$, and $\omega_j(s_i)$ are corresponding basis functions.

With the treatments above, $\ln\lambda_b$ can be represented by the values taken on a finite set of points. Therefore, the macroscopic design variables are collected in
\beq
\label{Eq:Macro_design_variables}
\textbf{F} = \left(\left(\textbf{F}^1\right)^{\top},\cdots,\left(\textbf{F}^M\right)^{\top}, \bar\theta\right)^{\top},
\eeq
where
\beq
\label{Eq:Macro_design_variables_per_curve}
\textbf{F}^i = \left(f^i_0,\cdots,f^i_{n_i-1}\right)^{\top}.
\eeq
\subsection{Sensitivity analysis}
\label{subsec:Sensitivity_analysis}
For any design variables, being either macroscopic or microscopic, say, $v \in\vD \cup \vF$, the adjoint method reads
\beq
\label{Eq:Derivative_compliance}
\pd{\mathcal{C}^{\text{H}}}{v} = -\int_{\Omega} \pd{\mathbb{C}_{ijkl}^{\text{H}}}{v} \pd{u_i^{\text{H}}}{x_j}\cdot \pd{u_k^{\text{H}}}{x_l}\dif \vx.
\eeq
Then we consider the derivatives of the homogenised elasticity tensor $\mathbb{C}^{\text{H}}$ with respect to design variables $v$. It is recalled from Eq.~\eqref{Eq:Conformal_cell_stiffness}~that $\mathbb{C}^{\text{H}}$ depends only on the rotation angle $\theta$ and the matrix cell $\hat{\mathbb{C}}^{\text{H}}$. Therefore, the evaluation of $\pd{\mathbb{C}_{ijkl}^{\text{H}}}{v}$ can be classified depending on the type of the design variables of interest, i.e.,
\beq
\label{Eq:Derivative_efficitive_property}
\pd{\mathbb{C}_{ijkl}^{\text{H}}}{v}=\left\{
\begin{aligned}
&R_{ip}R_{jq}R_{ks}R_{lt}\pd{\hat{\mathbb{C}}_{pqst}^{\text{H}}}{v}, &v \in \textbf{D} ,\\
&\frac{ \dif{\left(R_{ip}R_{jq}R_{ks}R_{lt}\right)} } { {\dif \theta} }\pd{\theta}{v}{\hat{\mathbb{C}}_{pqst}^{\text{H}}}, &v \in \textbf{F}.
\end{aligned}
\right.
\eeq
The first case of Eq.~\eqref{Eq:Derivative_efficitive_property} can be further calculated by {\cite{Xue_Smo2020}}
\beq
\label{Eq:Derivative_efficitive_property_repensent}
\pd{\hat{\mathbb{C}}_{pqst}^{\text{H}}}{v}=\int_{\Upsilon_{\text{p}}} \left( \delta_{pa}\delta_{qb}- \pd{\hat{\xi}_a^{pq}}{{\bar{Y}}_b} \right) \pd{\tilde{\mathbb{C}}_{abcd}^{\text{p}}}{v} \left(\delta_{sc}\delta_{td} -  \pd{   \hat{\xi}_c^{st}   }{\bar{Y}_d}\right) \dif \bar{\vY}, \quad  v \in \textbf{D}.
\eeq
Besides, since the rotation matrix $\textbf{R}=\left(R_{ij}\right)$ is an explicit function of $\theta$, $\frac{ \dif{\left(R_{ip}R_{jq}R_{ks}R_{lt}\right)} } { {\dif \theta} }$ should be expressed analytically.

We Incorporate Eqs.~\eqref{Eq:Normal_derivative}, ~\eqref{Eq:Arc_length_parameterisation_directional_derivative} and \eqref{Eq:Curve_lnlmd_parameterisation} into Eq.~\eqref{Eq:Green_function_neumann}, and obtain
\beq
\label{Eq:Derivative_Piecewise_Laplace_solution_Green_2}
\theta\left(\vx\right)=\sum_{i=0}^M \sum_{j=0}^{n_i}f^i_j\int_0^{l_i}\omega_j\left(s_i\right)G^{\textup{N}}\left(\vx,\vx'\right)\dif \Gamma_i + \bar\theta,\quad \vx' \in \Omega,
\eeq
where $G^{\textup{N}}\left(\vx,\vx'\right)$ and $\omega_j\left(s_i\right)$ are independent of the boundary values of $\ln\lambda_b$ given through $f^i_j$. Thus, $\theta\left(\vx\right)$ is actually a linear combination of the macroscopic design variables $f_j^i$. This means $\frac{\partial\theta}{\partial v}$ is effectively independent of the values of the macroscopic design variables. Consequently, $\frac{\partial\theta}{\partial v}$ in Eq.~\eqref{Eq:Derivative_Piecewise_Laplace_solution_Green_2} only needs to be solved once during the whole optimisation process, which further improves the computational efficiency of the algorithm.

\subsection{Optimisation procedure}
\label{subsec:Numerical_implementation}
\begin{figure}[!ht]
  \centering
    \includegraphics[width=.99\textwidth]{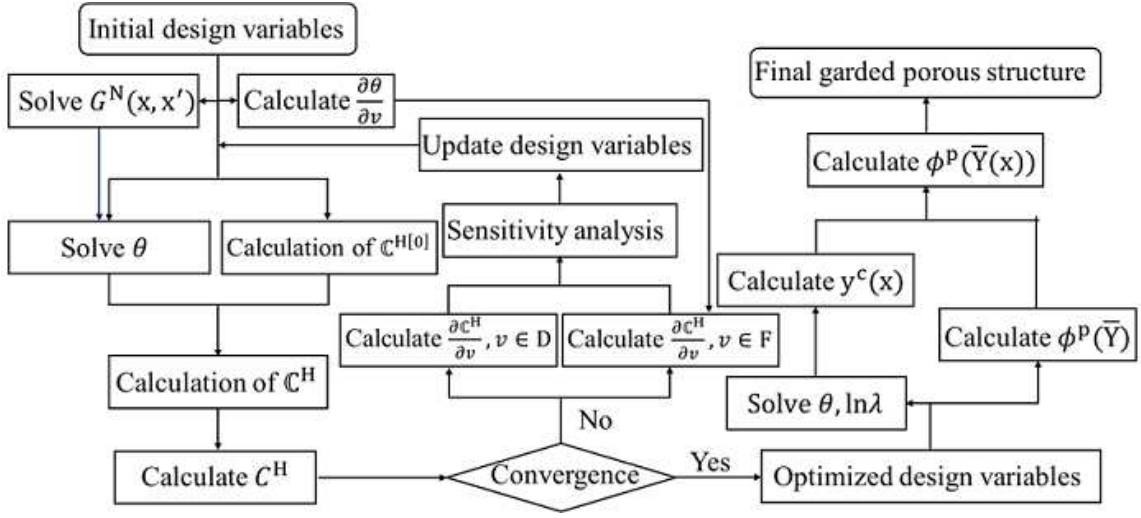}
    \caption{The flowchart of the AABH-plus-based optimisation framework combined with conformal mapping.}
    \label{fig:flowchart}
\end{figure}

Based on the derivation in the previous sections, the key steps to implement the optimisation algorithm presented here can be summerised as follows, and also by the flowchart in Fig.~\ref{fig:flowchart}.
\begin{enumerate}
    \item Initial stage
    \begin{enumerate}
        \item Initialise the microscopic and macroscopic design variables. Then calculate the Green's function $G^{\textup{N}}\left(\vx,\vx'\right)$ and $\frac{\partial \theta}{\partial v}, v\in \textbf{F}$ using Eq.~\eqref{Eq:Green_function_neumann} and Eqs.~\eqref{Eq:Derivative_Piecewise_Laplace_solution_Green_2}, respectively.
    \end{enumerate}
    \item Optimisation stage (Iteration until a converge criterion is met)
    \begin{enumerate}
      \item Solve for the first-order corrector $\hat{\bm{\xi}}$ and the homogenised elasticity tensor $\bm{\hat{\mathbb{C}}}$ through Eqs.~\eqref{Eq:Represent_cell_one_order_collector}--\eqref{Eq:Represent_stiffness_equivalent}. Determine the distribution of $\theta$ in $\Omega$ by solving Eq.~\eqref{Eq:Theta_solution}.
        \item Use Eq.~\eqref{Eq:Conformal_cell_stiffness}~to calculate the homogenised elasticity tensor $\bm{\mathbb{C}^{\text{H}}}$ in the macroscopic domain $\Omega$. Then analyse the resulting porous structure on the coarse grid to compute the equivalent system compliance $\mathcal{C}^{\text{H}}$.
        \item Evaluate the sensitivity quantities according to the combination of Eqs.~\eqref{Eq:Derivative_compliance}--\eqref{Eq:Derivative_efficitive_property_repensent} in Section \ref{subsec:Sensitivity_analysis}, where $\frac{\partial \theta}{\partial v},v \in \textbf{F}$ obtained in the initial stage can be used repeatedly without solving time and time again.
        \item Employ the methods of moving asymptotes (MMA) \cite{Svanberg_IJNME1987} as an optimisation tool to update the design variables.
    \end{enumerate}
    \item Representation of the optimised porous configuration
    \begin{enumerate}
  \item Determine the distribution of $\ln\lambda$ and $\theta$ by solving Eqs.~\eqref{Eq:Cauchy_Riemann_solution}--\eqref{Eq:Green_function} and then integrate Eq.~\eqref{Eq:Diff_y} to obtain the macroscopic mapping function $\vy^\text{c}\left(\vx\right)$ given by Eq.~\eqref{Eq:Integral_y}.
    \item Incorporate the macroscopic mapping function $\vy^\text{c}\left(\vx\right)$ and the TDF $\phi\left(\bar{\vY}\right)$ of the representing unit cell into Eq.~\eqref{Eq:Qusi_Peroidic_TDF}~to generate the optimised graded porous structure.
    \end{enumerate}
\end{enumerate}

Note that $\mathbb{C}^{\text{H}}$ is independent of $\ln\lambda$ at each point. Hence the system compliance can be obtained without solving \eqref{Eq:lnlmd_solution}. Therefore, $\ln\lambda$ only needs to be solved in the final stage when the optimised porous configuration is output.

\section{Numerical examples}
\label{sec: Numerical_examples}

In this section, several numerical examples based on the proposed optimisation framework are presented. Among them, we highlight the first example and compare it with the results by the projection method \cite{Groen_IJNME2018} to illustrate the effectiveness of the framework.

The examples discussed in this section are all two-dimensional in a plane stress state with unit thickness in the third dimension. Furthermore, we assume that the base materials are elastically isotropic and the parameters are taken as follows: Young's modulus $E =1$, and Poisson's ratio $\mu = 0.3$. Since the primary purpose of this paper is to investigate the advantages brought by conformal mapping, the matrix cell configuration is set fixed to be ``X'' shape with $\left|\Upsilon_{\text{p}}^{\text{s}}\right| = 0.3$. Besides, the scale separation parameter $\eps$ is kept to be 0.05, unless specified.

\subsection{Case with rectangle design domain}
\label{subsubsec:Rectangle_design_domain}
First, we consider a case where the design domain with length $L=2$ and width $W=1$, and it is discretised by a $100\times 50$ FEM mesh at the coarse-grained level. We choose 60 points on the boundary of the design domain and take the values of $\ln\lambda_b$ at these points as the design variables. The optimisation starts with a spatially periodic initial configuration. The upper and lower limits for $\ln\lambda_b$ are taken as $\pm\ln5$, which means the minimum and maximum of the scaling factor of the (unit) matrix cell are 0.2$\epsilon$ and 5$\epsilon$, respectively. Since the design domain is of rectangle shape, $\ln\lambda$ and $\theta$ can be solved for directly by the analytic formulae~\eqref{Eq:Rectangle_lnlmd_solution}-\eqref{Eq:Rectangle_theta_solution}.

\begin{figure}[!ht]
  \centering
  \includegraphics[width=.5\textwidth]{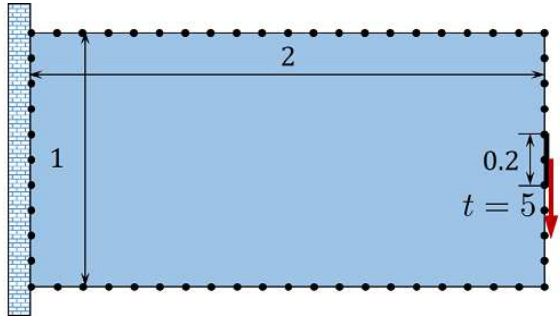}
  \caption{A cantilever beam subjected to a vertically distributed load over a small portion on the middle part of the right boundary.}\label{fig:Cantilever_beam}
\end{figure}

The loading scenario is selected as in Fig.~\ref{fig:Cantilever_beam}. A cantilever beam is fixed on its left side. A uniformly distributed load of $t=5$ is applied to the beam over a small portion on its right side. Such a set-up is to draw comparison with the benchmarked example obtained based on \cite{Groen_IJNME2018}. For comparative purposes, two ways of declaring the non-designable domains are investigated:
\begin{enumerate}[leftmargin={3.5em},label={\text{Case }\arabic*  }]
   \item A layer of solid materials with a thickness of 0.04 is deployed to the loading area.
   \item A layer of solid materials of a thickness of 0.02 is wrapped on the boundary of the design domain.
\end{enumerate}

\begin{figure}[!ht]
  \centering
    \subfigure[Case 1]{\includegraphics[width=.215\textwidth]{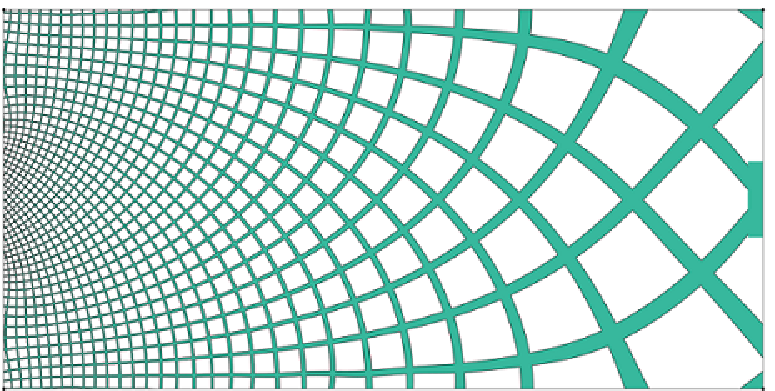}\label{fig:Optimised_cantilever_result_case1}}
    \quad
    \subfigure[Case 2]{\includegraphics[width=.215\textwidth]{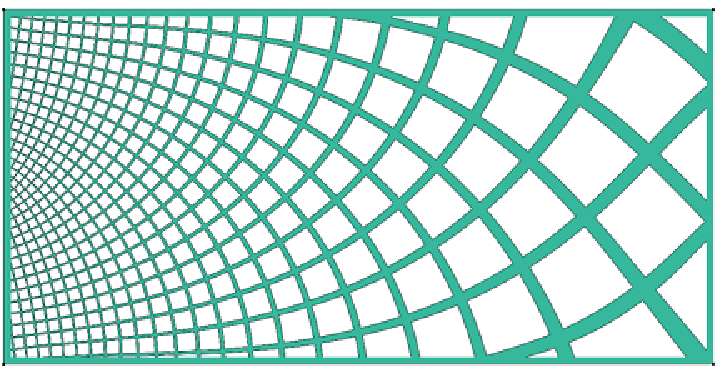}\label{fig:Optimised_cantilever_result_case2}}
  \quad
    \subfigure[Case 1 benchmark]{\includegraphics[width=.215\textwidth]{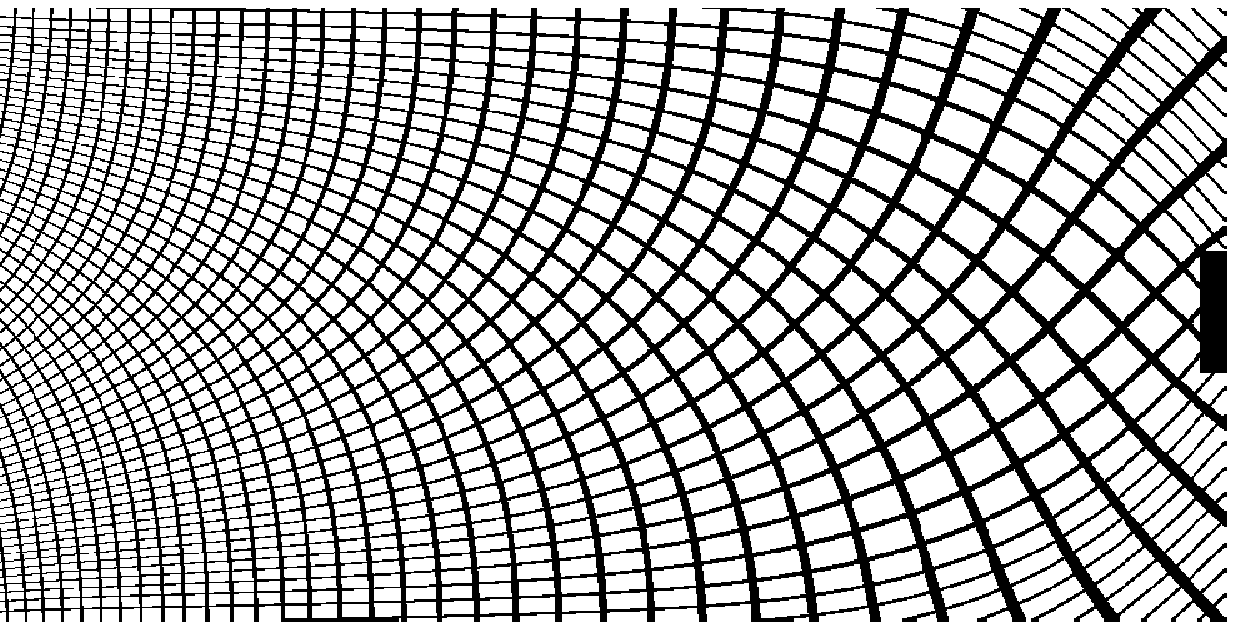} \label{fig:Optimised_cantilever_result_case1_projection}}
    \quad
    \subfigure[Case 2 benchmark]{\includegraphics[width=.215\textwidth]{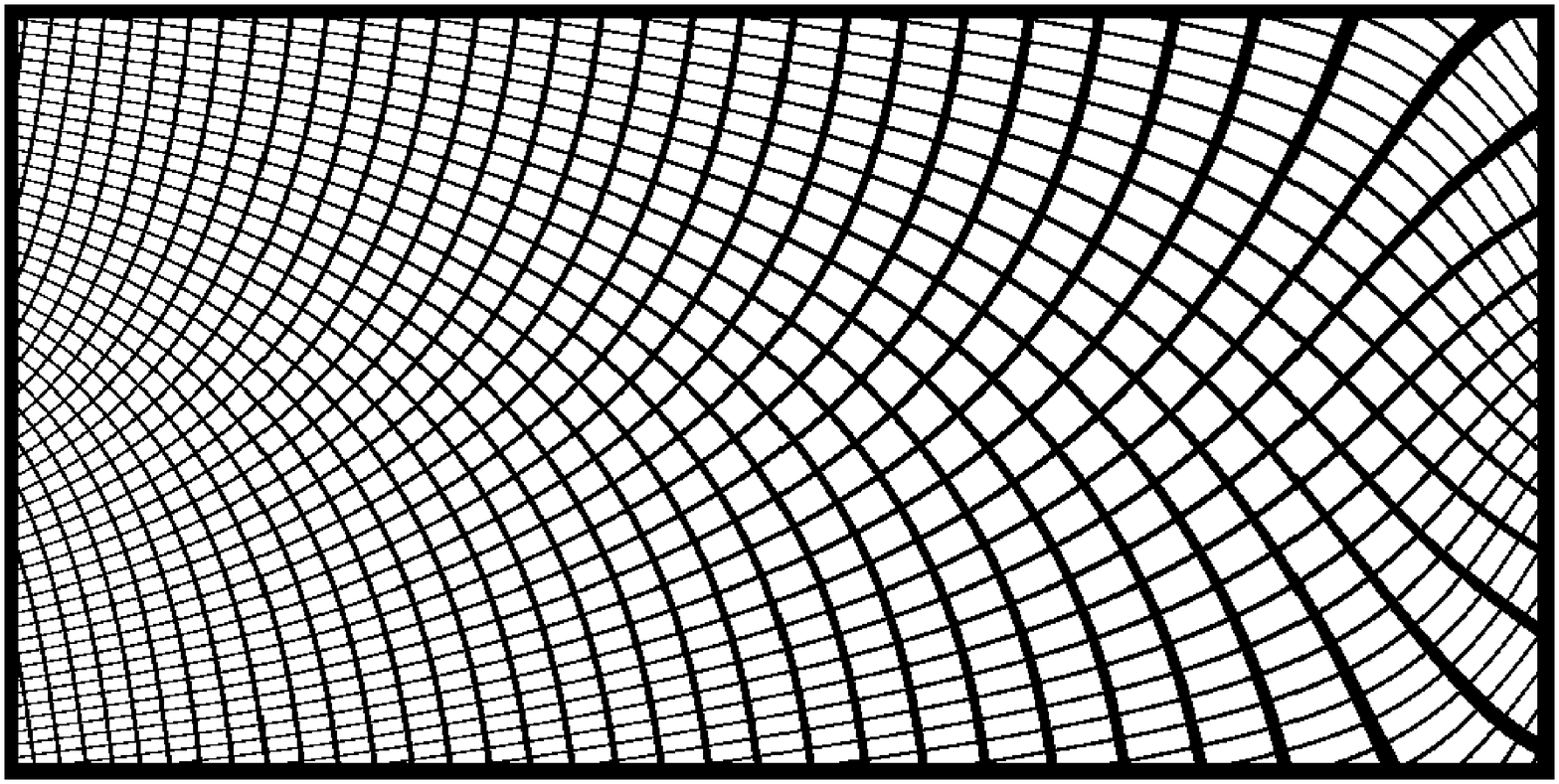} \label{fig:Optimised_cantilever_result_case2_projection}}
  \caption{The optimised results of the cantilever beam with case 1 and case 2 non-designable region: (a) and (b) based on the approach proposed in this study. (c) and (d) based on the projection approach \cite{Groen_IJNME2018}.}\label{fig:Optimised_cantilever_result}
\end{figure}

\begin{figure}[!ht]
  \centering
  \subfigure[Case 1]{\includegraphics[width=.46\textwidth]{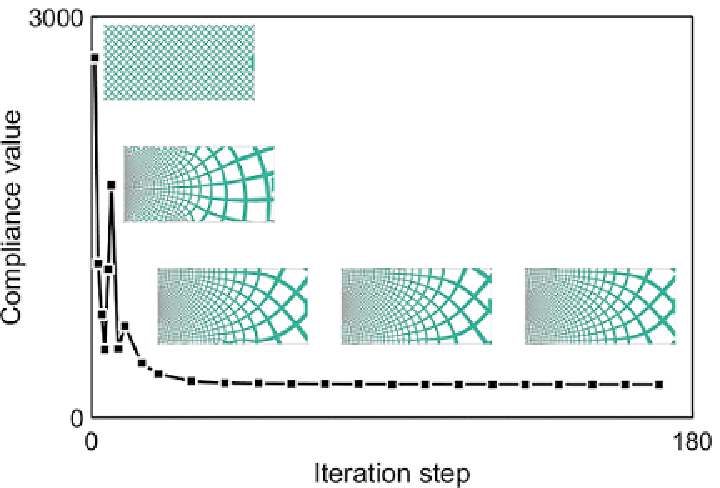}\label{fig:Case1_iteration_history}}
  \quad
  \subfigure[Case 2]{\includegraphics[width=.46\textwidth]{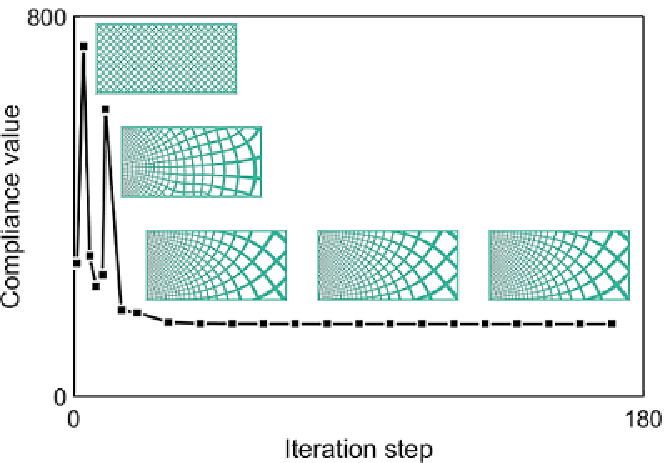}\label{fig:Case2_iteration_history}}
  \caption{The convergent histories of the compliance values calculated from the cantilever beam problems as shown in case 1 and case 2 in Fig.~\ref{fig:Optimised_cantilever_result}.}\label{fig:Optimised_cantilever_iteration_history}
\end{figure}

The optimised results for both cases are shown in Fig.~\ref{fig:Optimised_cantilever_result_case1} and \ref{fig:Optimised_cantilever_result_case2}, with the corresponding iteration histories shown in Fig.~\ref{fig:Optimised_cantilever_iteration_history}.  It can be read that the compliance values $\mathcal{C}^{\text{H}}$ in both cases become convergent after about 40 optimisation steps and the final results are 242.67 and 151.31, respectively.

To check the accuracy of the multiscale method in this framework, we analyse the optimisation results on a fine mesh with $2000 \times 1000$ elements. The obtained compliance values from fine-scale simulation denoted by $\mathcal{C}^{\text{FS}}$ are very close to the corresponding compliance values based on the AABH plus approach, as shown in Table.~\ref{Tap:Comparision_compliances}. Furthermore, we compare the optimisation results with their counterparts obtained by using the projection method \cite{Groen_IJNME2018} as shown in Fig.~\ref{fig:Optimised_cantilever_result_case1_projection} and~\ref{fig:Optimised_cantilever_result_case2_projection}. The corresponding compliance values of both homogenisation-based and fine-scale analysis are comparatively shown in Table.~\ref{Tap:Comparision_compliances}. It can be read from the present method that the results are highly consistent with that of the projection method \cite{Groen_IJNME2018}, in the aspects of both the configuration outlook and the system compliance.
\begin{table}[htbp]
  \centering
  \caption{Comparison of the compliances of optimised results based on the AABH-plus-based framework combined with conformal mapping and their counterparts based on the projection approach \cite{Groen_IJNME2018}.}\label{Tap:Comparision_compliances}
  \vspace{0.2cm}
  \begin{tabular}{ccccc}
  \toprule
       & \multicolumn{2}{c}{AABH plus}   & \multicolumn{2}{c}{Projection}\\
  \cmidrule{2-5}
        & $\mathcal{C}^\text{H}$ & $\mathcal{C}^\text{FS}$   & $\mathcal{C}^\text{H}$  & $\mathcal{C}^\text{FS}$  \\
  \midrule
  Case 1     & 242.67    & 249.94    &228.04   &242.26\\
  Case 2     & 151.31    & 155.66    &152.21   &158.11\\
  \bottomrule
  \end{tabular}
  \end{table}

Besides, we also examine the feature size control of this example with case 2. To this end, a series of optimised results with different combinations of $\lambda_{b\max}$ and $\lambda_{b\min}$ are obtained. All the final configurations and corresponding compliance values, including $\mathcal{C}^\text{H}$ and $\mathcal{C}^\text{FS}$, are listed in Table.~\ref{Tab:Case2_result_diffierent_lambda}. From this table, it can be read that the minimum characteristic size of the porous structure increases with an increase in $\lambda_{b\min}$, but the deviation of $\mathcal{C}^\text{H}$ from $\mathcal{C}^\text{FS}$ is not affected much. Besides, the computational time $T$ in the MATLAB environment with a serial programmer on a desktop platform for each situation is also shown in Table.~\ref{Tab:Case2_result_diffierent_lambda}. It can be seen that the whole optimisation process roughly takes about 100s, which makes a big improvement compared with the preceding study using a zoning scheme, which delivers a better performance with computational parallelism \cite{Xue_Smo2020}.
\begin{table}[!ht]
  \centering
  \caption{
  The optimised results with different combinations of $\lambda_{b\min}$ and $\lambda_{b\max}$. Here $\lambda_{b\min}$ determines the minimum size of the microstructural members.
  }\label{Tab:Case2_result_diffierent_lambda}
  \vspace{0.2cm}
  \rowcolors{2}{black!10}{black!5}
  \renewcommand{\arraystretch}{1.5}
  \begin{tabular}{m{1.2cm}<{\centering}m{1.2cm}<{\centering}m{1.2cm}<{\centering}m{1.2cm} <{\centering}m{1.2cm}<{\centering}m{1.2cm}<{\centering}m{1.3cm}<{\centering}m{3.0cm} <{\centering}}
\firsthline
  \rowcolor{black!20}
  $\eps$  &$\lambda_{b\min}$ & $\lambda_{b\max}$& $V_f$ &  $\mathcal{C}^{\text{H}}$ & $\mathcal{C}^{\text{FS}}$ &$T$ & \\
 \hline
    0.05 &1/4 &4 &0.337 & 151.60 &157.21 & 96.69s &
    \begin{minipage}[c]{0.18\textwidth}
      \centering
      \includegraphics[width=\textwidth]{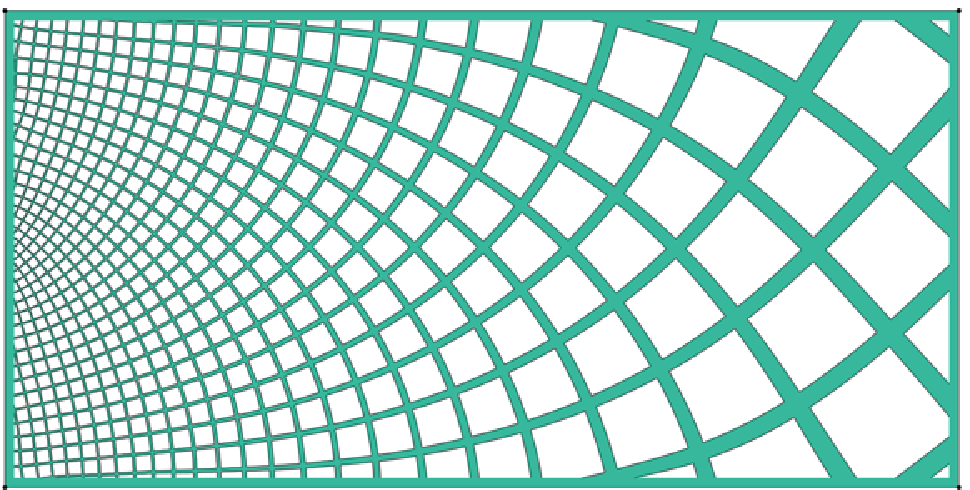}
    \end{minipage}\\
    0.05 &1/3 &4 & 0.337    & 151.59 &157.57 & 99.29s &
    \begin{minipage}[c]{0.18\textwidth}
      \centering
      \includegraphics[width=\textwidth]{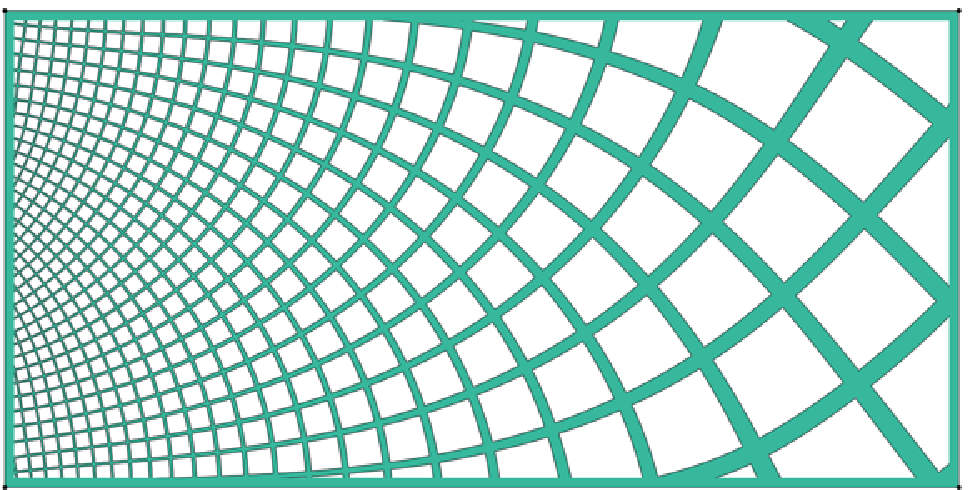}
    \end{minipage}\\
    0.05 &1/2 &4 & 0.336 & 155.43 &169.14 &100.26s &
    \begin{minipage}[c]{0.18\textwidth}
      \centering
      \includegraphics[width=\textwidth]{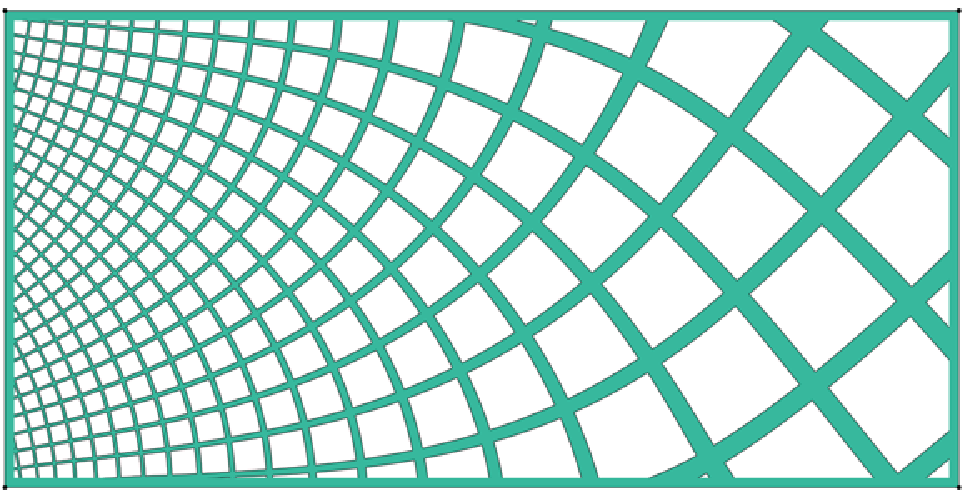}
    \end{minipage}\\
    0.05 &2/3 &4 & 0.336 & 162.73 &171.11 &99.87s &
    \begin{minipage}[c]{0.18\textwidth}
      \centering
      \includegraphics[width=\textwidth]{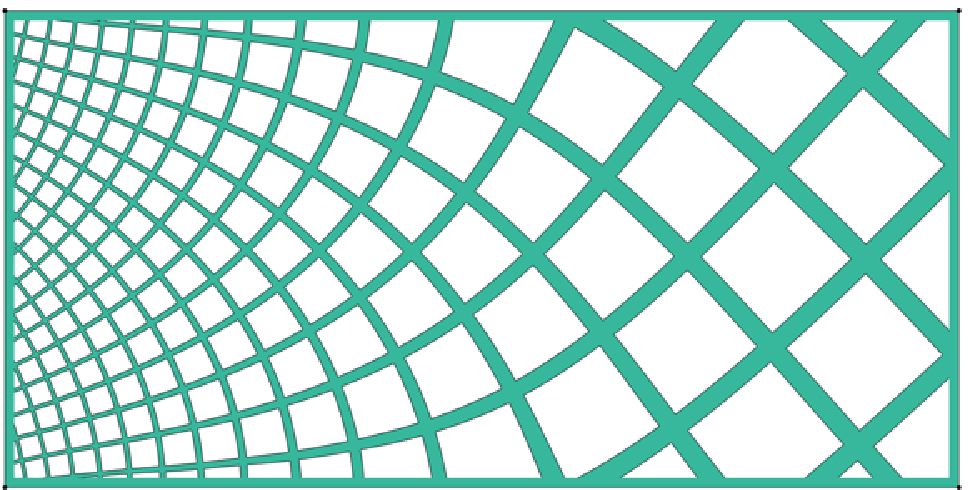}
    \end{minipage}\\
    0.05 &1   &4 & 0.337 & 179.15 &188.87 &101.61s &
    \begin{minipage}[c]{0.18\textwidth}
      \centering
      \includegraphics[width=\textwidth]{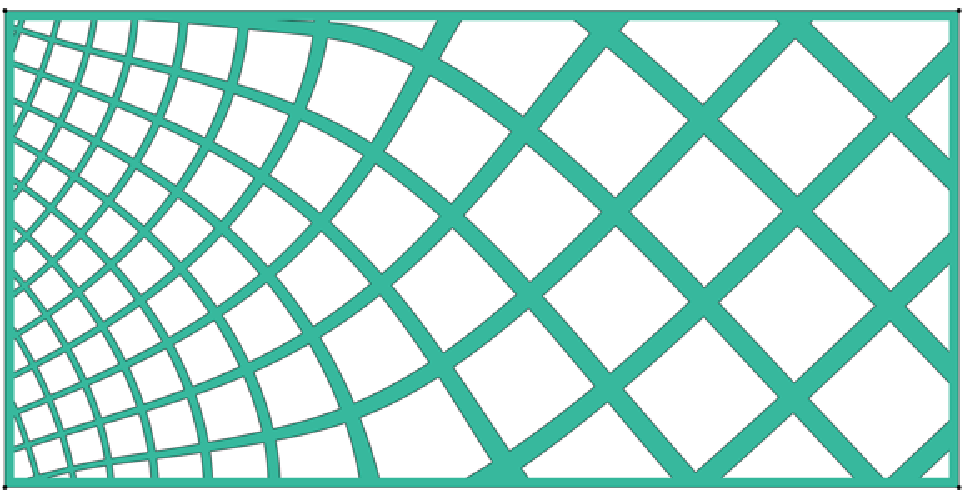}
    \end{minipage}\\
\hline
  \end{tabular}
\end{table}

\begin{figure}[!ht]
  \centering
  \includegraphics[width=.5\textwidth]{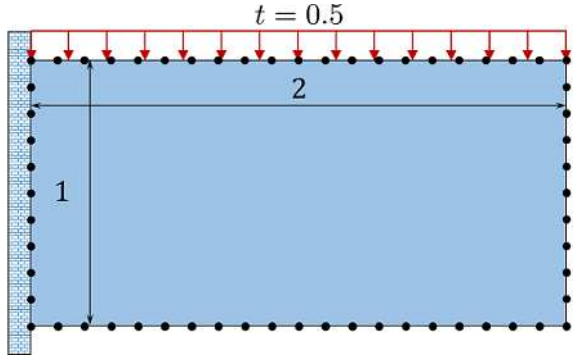}
  \caption{A cantilever beam subjected to a vertically distributed load on the top side.\label{fig:Cantilever_beam_up}}
\end{figure}

\begin{figure}[!ht]
  \centering
  \subfigure[]{\includegraphics[width=.48\textwidth]{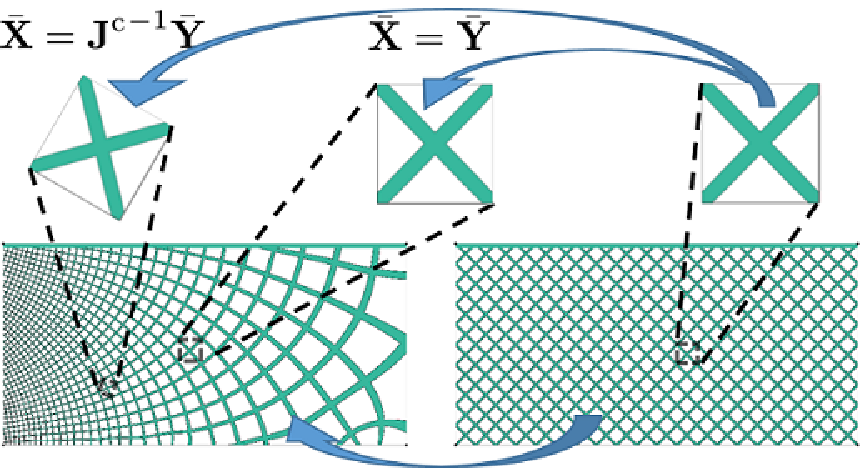} \label{fig:Cantilever_beam_up_result}}
  \quad
  \subfigure[]{\includegraphics[width=.46\textwidth]{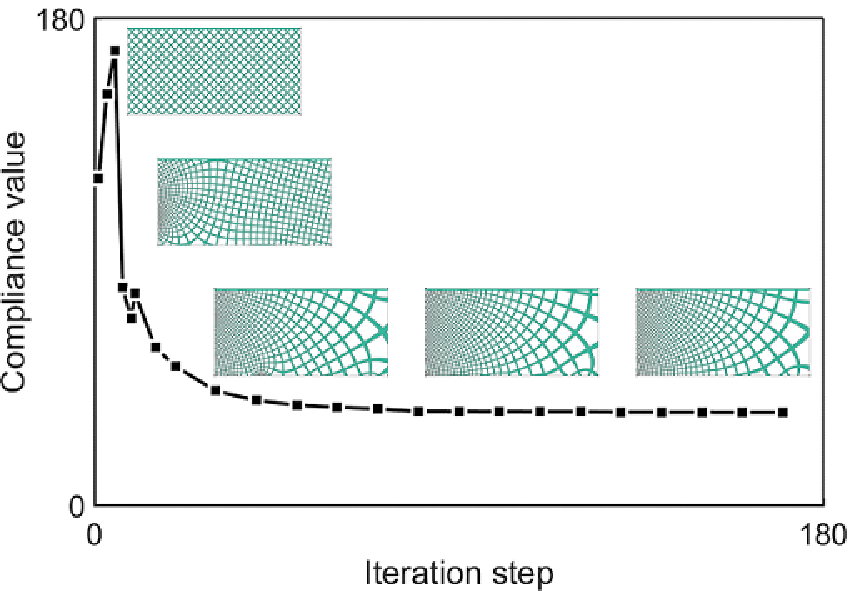} \label{fig:Cantilever_beam_up_iteration_history}}
  \caption{A cantilever beam with a distributed load on the top side: (a) Optimised result. (b) Compliance values convergence history. \label{fig:Optimised_cantilever_beam_up_result}}
\end{figure}

\begin{figure}[!ht]
  \centering
  \includegraphics[width=.5\textwidth]{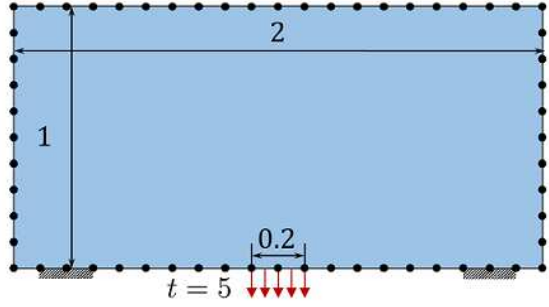}
  \caption{A bridge subjected to a vertically distributed load in the middle part of the bottom boundary.\label{fig:Bridge}}
\end{figure}

\begin{figure}[!ht]
  \centering
  \subfigure[]{\includegraphics[width=.49\textwidth]{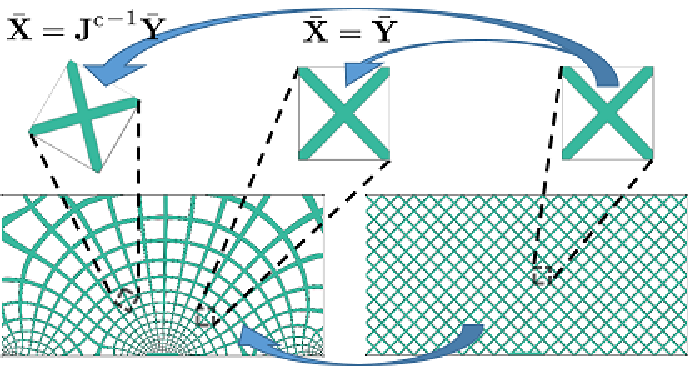}\label{fig:Bridge_result}}
  \quad
  \subfigure[]{\includegraphics[width=.46\textwidth]{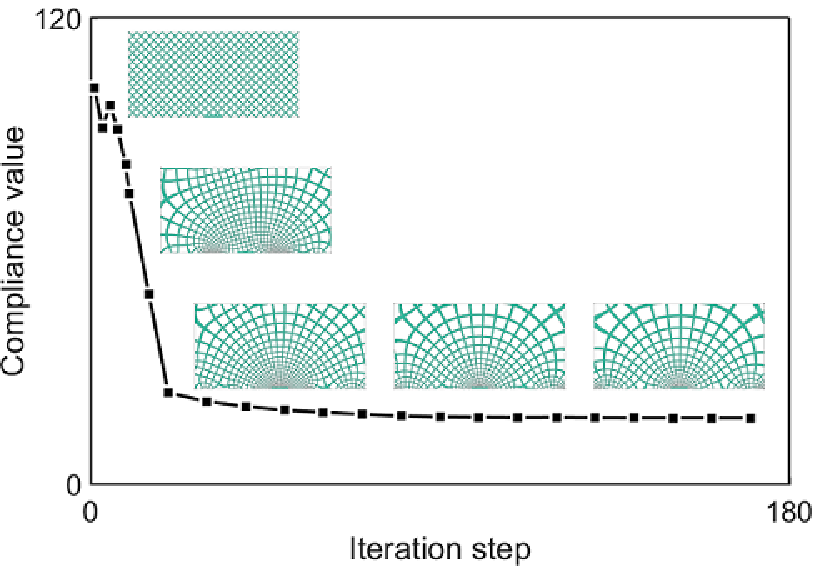} \label{fig:Bridge_iteration_history}}
  \caption{A bridge problem: (a) Optimised result. (b) Compliance values convergence history. \label{fig:Optimised_Bridge_result}}
\end{figure}

Other compliance minimisation problems are also considered: A cantilever beam subjected to a vertically distributed load $t = 0.5$ on the top side, as shown in Fig.~\ref{fig:Cantilever_beam_up}. The corresponding optimised configuration and the compliance convergence curve are shown in Fig.~\ref{fig:Optimised_cantilever_beam_up_result}. We also consider the case of a bridge with a vertically distributed load $t= 5$ in the middle part of the bottom boundary as shown in Fig.~\ref{fig:Bridge}. The corresponding result and the iteration history are shown in Fig.~\ref{fig:Optimised_Bridge_result}. In both situations, sensible results are obtained.

In summary, for the case with rectangular design domains, since the solution to Laplace's equation can be represented analytically, the partial differential equations that need to be solved are just Eqs.~\eqref{Eq:Equilibrium_equation} and \eqref{Eq:Represent_cell_one_order_collector} during each iteration step. Especially when the configuration of the representing unit cell remains unchanged, the first-order corrector $\hat{\bm{\xi}}$ in Eq.~\eqref{Eq:Represent_cell_one_order_collector} can be solved only once before optimisation. Consequently, only the homogenised equilibrium equation Eq.~\eqref{Eq:Equilibrium_equation} needs to be solved on the coarse mesh in the optimisation process. This maximally brings down the computational cost of the present method, making its computation efficiency be levelled with general scale-separation schemes, while the model accuracy is greatly improved.

\subsection{Design domain of general shape}
\label{subsubsec:General_design_domain}
\begin{figure}[!ht]
  \centering
  \includegraphics[width=.5\textwidth]{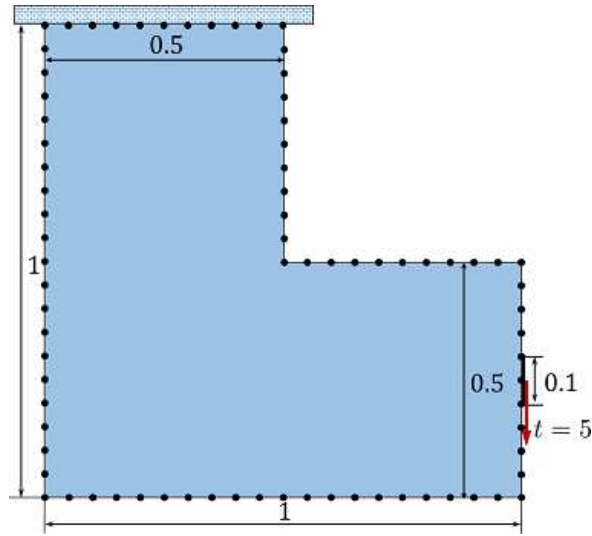}
  \caption{An L-beam subjected to a vertically distributed load in the middle part of the right boundary.\label{fig:L_beam}}
\end{figure}

\begin{figure}[!ht]
  \centering
  \subfigure[]{\includegraphics[width=.45\textwidth]{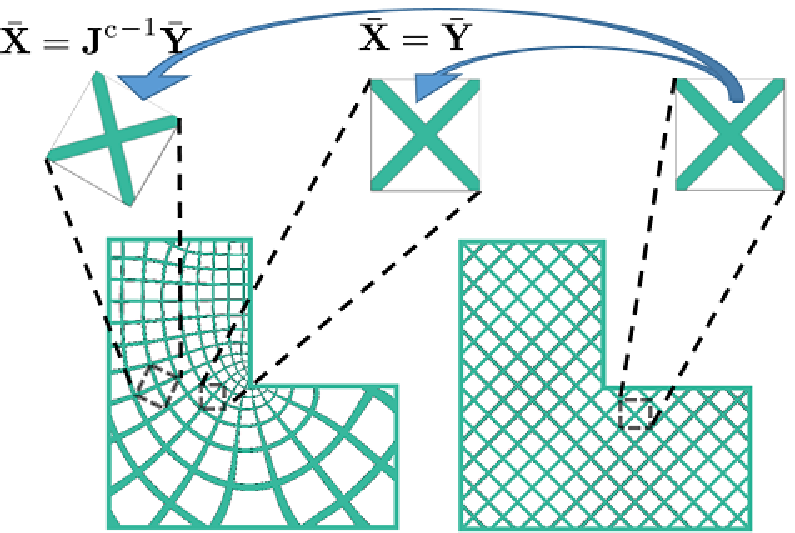} \label{fig:L_beam__result}}
  \quad
  \subfigure[]{\includegraphics[width=.47\textwidth]{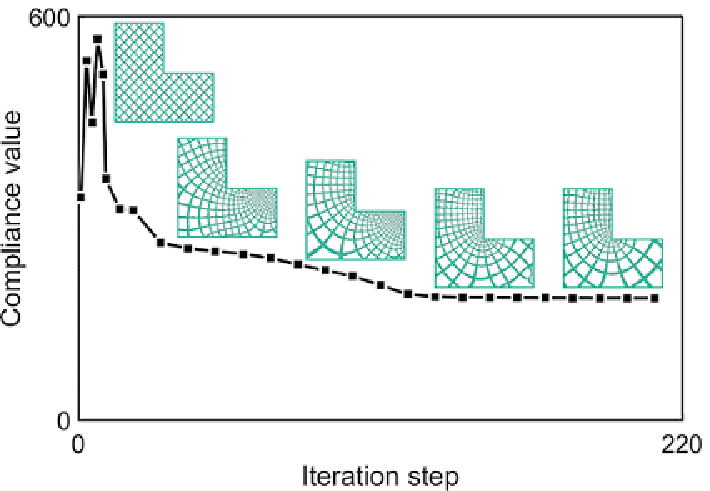} \label{fig:L_beam_iteration_history}}
  \caption{An L-beam problem: (a) The optimised result shown on the left and the initial configuration shown on the right. (b) Compliance values convergence history. \label{fig:Optimised_L_beam__result}}
\end{figure}

Now, we consider the case, where the shape of the design domain is more general. To this end, we consider an L-beam with a vertically distributed load $t=5$ applied at the middle point of its right side, as shown in Fig.~\ref{fig:L_beam}. This specimen is fixed on the top side. We choose 80 points on the boundary, as shown in Fig.~\ref{fig:L_beam}, as the nodes at which the design variables get evaluated. The corresponding bounds for $\ln\lambda$ are set to be $\pm \ln5$. The distributions of $\ln\lambda$ and $\theta$ field in the domain are solved for by using Eqs.~\eqref{Eq:Cauchy_Riemann_solution}-\eqref{Eq:Green_function}. The optimisation result is shown in Fig.~\ref{fig:Optimised_L_beam__result}. Again, sensible results can be observed. Note that large stress concentration is expected at the upper right corner of the L-shape beam, and the microstructures there rotate fast so as to levy the highly stressed state there.

\section{Discussion}
\label{sec: Discussion}
\subsection{Summary of the method}
\label{subsec: Summary}
In this article, with the use of the holomorphic functions as the macroscopic mapping function, the concept of conformal mapping is employed to improve the AABH-plus-based optimisation framework where the representation and optimisation of graded orthotropic porous structures are effectively facilitated. The presented approach not only inherits the merits of both conformal mapping and the AABH-plus-based framework such as angle-preserving, smoothness, high gradual flexibility, etc but also exhibits several advantage features summarised as follows. Firstly, compared with general cases where cell problems have to be resolved point-wisely, here one just needs to solve for the elasticity tensors corresponding to the matrix cells, and the computational efficiency thus gets significantly improved. For some examples, the whole optimisation process just costs about 100s on a normal desktop computer. Secondly, due to the maximum principle of Laplace's equation, the minimal characteristic length of the porous configuration of interest can be monitored by directly limiting the upper and lower bound of the macroscopic design variables, that is, the logarithm values of the scaling factor at a collection of selected points on the boundary. Thirdly, the rotation angle field in this framework is naturally continuous without extra process in case of angle jumps of $\pi$. Through the digital representation of design variables and corresponding sensitivity analysis, the numerical implementation framework of the presented approach is established. Several numerical examples are considered and some of them are compared with benchmarks based on other recently developed methods.

\subsection{Comparative discussion with other methods in the frontier}
\label{subsec: Comparative}
The similarities and differences between the proposed method here and other novel methods for the design of porous structures \cite[e.g.][]{Pantz_SIAM_JCO2008, Vogiatzis_CMME2018, Groen_IJNME2018, Allaire_CMA2019, Zhu_JMPS2019, Jiang_FME2019,Xie_CAD2020} are worth being discussed. It has been summarised \cite{Xue_CMME2020} that all the existing approaches for representing porous configurations can be effectively reformulated through the mapping approach, as shown in Fig.~\ref{fig:GMs_generation}. This enables one to examine the aforementioned studies, as well as the present one, from a unified viewpoint. In the studies, where conformal mapping is somehow involved \cite{Vogiatzis_CMME2018, Jiang_FME2019}, the effects due to cell rotation angle on structural were barely considered. As for the projection-based method \cite[e.g.][]{Pantz_SIAM_JCO2008, Groen_IJNME2018, Allaire_CMA2019,Xie_CAD2020}, the constituting cell rotation angles at different grids are all set to be mutually independent design variables. Then the projection method is implemented to produce the final configuration whose constituting cell direction is as consistent as possible with the orientation field obtained during the optimisation stage. Here two issues are worth being addressed. Firstly, since the angles at two neighbouring pixels are irrelevant to each other, there may be an angle jump of $\pi$ between them given the rotational
symmetry shown by the homogenised elasticity tensor. Secondly, the optimised cell rotation angle field, i.e, the gradient direction of the mapping function, usually does not take a total differential form. Therefore, the real angle field of the final configuration generated by the projection method in the post-process stage becomes an approximation to the designed mapping function. In contrast, in the proposed framework the orientation field determined by C-R equations automatically meet the exact differential condition without post-process. Hence a smooth angle field without jumps of angle $\pi$ between neighbouring pixels is naturally ensured. Here, not the rotation angle, but the logarithm value of the scaling factor of the constituting cell is chosen to be the design controller with the conditions for total differential being naturally met, the TDF for the optimal porous configuration is at hand. Such a strategy combined with the maximum principle of Laplace's equation realised the minimal feature size control in an extraordinarily simple way.
\subsection{Implication on further studies}
\label{subsec:future}
Further studies are expected along the following two directions. Firstly, the presented approach needs to be generalised for wider coverage. For instance, the work by Xue et al. \cite{Xue_CMME2020} can be incorporated into this framework to generate graded porous configurations, where the volume fraction of constituting cell is allowed to gradually vary. Furthermore, the disappearance mechanism of the cells can be introduced to obtain optimised results possessing macroscopic topologies, realising concurrent optimisation on different scales. Moreover, graded porous structures containing a set of desired seminal cells also could be achieved and optimised. In addition, this method may also be extended in other fields such as acoustics and thermal. Secondly, the presented method also needs to be deployed for three-dimensional cases. In this scenario, it can be seen that conformal mapping is fully restricted to a composition of isometries, dilatations or inversions, at most one of each, based on Liouville's theorem \cite{Carmo_1992}. This means using conformal mapping as the macroscopic mapping will severely limit the design freedom of graded porous structures for three-dimensional cases. Therefore, new types of mapping functions, such as quasi-conformal mapping, are still expected. They should enjoy a higher degree of design freedom than conformal mapping, while serious cell distortion is still avoided.

\section*{Acknowledgement}
We thank Ole Sigmund and Jeroen Groen for providing the comparative simulation examples shown in Fig.~\ref{fig:Optimised_cantilever_result}(c) and (d). The financial supports from National Key Research and Development Plan (2020YFB1709400) from the Ministry of Science and Technology of the People's Republic of China, the National Natural Science Foundation of China (11772076, 11675161, 11732004, 11821202) are gratefully acknowledged.

\section*{References}

\bibliography{myreference}
  \newcounter{appendixcounter}
  \setcounter{appendixcounter}{0}
  \renewcommand\theappendixcounter{\Alph{appendixcounter}}
\refstepcounter{appendixcounter}
\section*{Appendix \theappendixcounter}
\label{Appendix:Cauchy_Riemann}
In this appendix, that $\ln\lambda$ and $\theta$ satisfy Eq.~\eqref{Eq:Cauchy_Riemann} is proven. First, expanding Eq.~\eqref{Eq:Integrable}~gives
\begin{subequations}\label{Eq:Integrable_expansion}
    \begin{alignat}{2}
        \label{Eq:Integrable_expansion_1}
        &\left(\frac{1}{\lambda^2}\pd{\lambda}{x_1}-\frac{1}{\lambda}\pd{\theta}{x_2}\right) \sin\theta-\left(\frac{1}{\lambda^2}\pd{\lambda}{x_2}+\frac{1}{\lambda} \pd{\theta}{x_1}\right) \cos\theta=0,\\
        \label{Eq:Integrable_expansion_2}
        &\left(\frac{1}{\lambda^2}\pd{\lambda}{x_1}-\frac{1}{\lambda}\pd{\theta}{x_2}\right) \cos\theta+\left(\frac{1}{\lambda^2}\pd{\lambda}{x_2}+\frac{1}{\lambda} \pd{\theta}{x_1}\right)\sin\theta=0.
    \end{alignat}
\end{subequations}
Then, re-writing Eq.~\eqref{Eq:Integrable_expansion} in matrix form gives
\beq
\label{Eq:Integrable_expansion_matrix_form}
\frac1{\lambda}
\begin{pmatrix}
    \sin\theta &  -\cos\theta\\
    \cos\theta &    \sin\theta
\end{pmatrix}
\begin{pmatrix}
    \frac{1}{\lambda}\frac{\partial \lambda}{\partial x_1}-\frac{\partial \theta}{\partial x_2}\\
    \frac{1}{\lambda}\frac{\partial{\lambda}}{\partial x_2}+\frac{\partial \theta}{\partial x_1}
\end{pmatrix}=
\begin{pmatrix}
    0\\
    0
\end{pmatrix}.
\eeq

Obviously, $\frac{1}{\lambda}\begin{pmatrix}
    \sin\theta &  -\cos\theta\\
    \cos\theta &    \sin\theta
\end{pmatrix}$ is invertible, thus
\begin{subequations}
\label{Eq:Cauchy_Riemann_medium}
\begin{alignat}{2}
    \label{Eq:Cauchy_Riemann_medium_1}
    &\frac{1}{\lambda}\frac{\partial \lambda}{\partial x_1}-\frac{\partial \theta}{\partial x_2} = 0,\\
    \label{Eq:Cauchy_Riemann_medium_2}
    &\frac{1}{\lambda}\frac{\partial{\lambda}}{\partial x_2}+\frac{\partial \theta}{\partial x_1} = 0.
\end{alignat}
\end{subequations}
By substituting $\frac1{\lambda}\frac{\partial \lambda}{\partial x_i} = \frac{\partial \ln\lambda}{\partial x_i},i = 1,2$ into Eq.~\eqref{Eq:Cauchy_Riemann_medium}, Eq.~\eqref{Eq:Cauchy_Riemann} is proven.

\refstepcounter{appendixcounter}
\section*{Appendix \theappendixcounter}
\label{Appendix:One_order_corrector}
This appendix shows $\bm{\xi}$ and $\mathbb{C}^{\text{H}}$ satisfy Eqs.~\eqref{Eq:Conformal_cell_one_order_collector} and \eqref{Eq:Conformal_cell_stiffness}, respectively. Since the base materials are isotropic, whose elasticity tensor has the same components in any reference frame,
\beq
\label{Eq:Stiffness_tensor_rotation}
\tilde{\mathbb{C}}^{\text{p}}_{ijkl}=R_{ia}R_{jb}R_{kc}R_{ld} \tilde{\mathbb{C}}^{\text{p}}_{abcd}.
\eeq
Replacing the indices $k,l,c$ and $d$ with $s,t,u$ and $v$, Eq.~\eqref{Eq:Stiffness_tensor_rotation}~becomes
\beq
\label{Eq:Stiffness_tensor_rotation_2}
\tilde{\mathbb{C}}^{\text{p}}_{ijst}=R_{ia}R_{jb}R_{su}R_{tv} \tilde{\mathbb{C}}^{\text{p}}_{abuv}.
\eeq

It is recalled from Eq.~\eqref{Eq:J_R} that the components of the Jacobian $\vJ^{\text{c}}$ and the rotation matrix $\textbf{R}$ satisfy
\beq
\label{Eq:Rotation_matrix_components}
{J}_{ij}=\frac{1}{\lambda}{R}_{ji}.
\eeq

Incorporating Eqs.~\eqref{Eq:Stiffness_tensor_rotation}-\eqref{Eq:Rotation_matrix_components} into Eq.~\eqref{Eq:One_order_collector} gives
\beq
\label{Eq:One_order_corrector_midium_1}
\frac{1}{\lambda}\pd{}{\bar{{Y}}_m} \left(R_{jm} R_{ln}R_{ia}R_{jb}R_{kc}R_{ld}\tilde{\mathbb{C}}^{\text{p}}_{abcd} \pd{{\xi_k^{st}}}{\bar{{Y}}_n}\right)=R_{jm}R_{ia}R_{jb}R_{su}R_{tv}
\pd{\tilde{\mathbb{C}}^{\text{p}}_{abuv} }{\bar{Y}_m}.
\eeq

Since $\textbf{R}$ is an orthogonal matrix,
\beq
\label{Eq:Orthogonal_matrix_production}
R_{jb}R_{jm}={\delta}_{bm},\quad R_{ld}R_{ln}={\delta}_{dn},
\eeq
where $\bm{\delta} = \left(\delta_{ij}\right)$ is the Kronecker delta. With the use of Eq.~\eqref{Eq:Orthogonal_matrix_production},
Eq.~\eqref{Eq:One_order_corrector_midium_1} reduces to
\beq
\label{Eq:One_order_corrector_midium_2}
\frac{1}{\lambda}\pd{}{\bar{{Y}}_b} \left(\tilde{\mathbb{C}}^{\text{p}}_{abcd}R_{ia}R_{kc}\pd{\xi_k^{st}}{\bar{{Y}}_d}\right)= R_{ia}R_{su}R_{tv}\pd{\tilde{\mathbb{C}}^{\text{p}}_{abuv} }{\bar{Y}_b}.
\eeq

Multiplying both sides of Eq.~\eqref{Eq:One_order_corrector_midium_2}~by $R_{io}R_{se}R_{tf}$ and re-using Eq.~\eqref{Eq:Orthogonal_matrix_production}~gives
\beq
\label{Eq:One_order_corrector_midium_3}
\frac{1}{\lambda}\pd{}{\bar{{Y}}_b} \left(\tilde{\mathbb{C}}^{\text{p}}_{obcd}R_{kc}R_{se}R_{tf}\pd{\xi_k^{st}}{\bar{{Y}}_d}\right)= \pd{\tilde{\mathbb{C}}^{\text{p}}_{obef}}{\bar{Y}_b}.
\eeq
Using indices $i,j,k,l,u,v,s,t$ and $c$ instead of $o,b,c,d,s,t,e,f$ and $k$, Eq.~\eqref{Eq:One_order_corrector_midium_3} is re-written by
\beq
\label{Eq:One_order_corrector_midium_4}
\frac{1}{\lambda}\pd{}{\bar{{Y}}_j} \left(\tilde{\mathbb{C}}^{\text{p}}_{ijkl}R_{ck}R_{us}R_{vt}\pd{\xi_c^{uv}}{\bar{{Y}}_l}\right)= \pd{\tilde{\mathbb{C}}^{\text{p}}_{ijst} }{\bar{Y}_j}.
\eeq

Comparing Eq.~\eqref{Eq:One_order_corrector_midium_4} to Eq.~\eqref{Eq:Represent_cell_one_order_collector}, gives
\beq
\label{Eq:inverse_one_order_corrector}
\hat{\xi}_k^{st}=\frac{1}{\la}\xi_c^{uv}R_{ck}R_{us}R_{vt}.
\eeq
Multiplying both sides of Eq.~\eqref{Eq:inverse_one_order_corrector} by $\lambda R_{ok}R_{ps}R_{qt}$ and using Eq.~\eqref{Eq:Orthogonal_matrix_production} again gives
\beq
\label{Eq:One_order_corrector_midium_5}
\xi_o^{pq} = \lambda R_{ps}R_{qt} R_{ok}\hat{\xi}_k^{st}.
\eeq
Letting the indices $u$, $v$ and $w$ in place of $p$, $q$ and $o$ as free indices in Eq.~\eqref{Eq:One_order_corrector_midium_5}, Eq.~\eqref{Eq:Conformal_cell_one_order_collector} is proven.

Replacing the indices $u,v,w,s,t$ and $k$ in Eq.~\eqref{Eq:Conformal_cell_one_order_collector} with $k,l,s,c,d$ and $e$ gives
\beq
\label{Eq:Conformal_cell_one_order_collector_2}
{\xi_s^{kl}}={\lambda}R_{se}R_{kc}R_{ld}{\hat{\xi}_e^{cd}}.
\eeq
Incorporating Eq.~\eqref{Eq:Stiffness_tensor_rotation}--\eqref{Eq:Rotation_matrix_components} and~\eqref{Eq:Conformal_cell_one_order_collector_2}~into Eq.~\eqref{Eq:stiffness_equivalent} gives
\beq
\label{Eq:stiffness_equivalent_2}
{\mathbb{C}_{ijkl}^{\text{H}}}
= R_{ia}R_{jb}R_{kc}R_{ld}{{\mathbb{C}}}_{abcd}\cdot |\Upsilon_{\text{p}}^\text{s}| - R_{tn}R_{se}R_{kc}R_{ld}R_{ia}R_{jb}R_{su}R_{tv}\int_{\Upsilon_{\text{p}}^{\text{s}}}
{{\mathbb{C}}}_{abuv}  \pd{\hat{\xi}_e^{cd}}{\bar{Y}_n} \dif \bar{\vY}.
\eeq
According to Eq.~\eqref{Eq:Orthogonal_matrix_production} gives
\beq
\label{Eq:Orthogonal_matrix_production_2}
R_{se}R_{su}={\delta}_{eu},\quad R_{tn}R_{tv}={\delta}_{nv}.
\eeq
By substituting Eq.~\eqref{Eq:Orthogonal_matrix_production_2} into Eq.~\eqref{Eq:stiffness_equivalent_2}, $\mathbb{C}_{ijkl}^{\text{H}}$ is re-written by
\beq
\label{Eq:Orthogonal_matrix_production_3}
\mathbb{C}_{ijkl}^{\text{H}}
= R_{ia}R_{jb}R_{kc}R_{ld}\left({{\mathbb{C}}}_{abcd}\cdot |\Upsilon_{\text{p}}^\text{s}| - \int_{\Upsilon_{\text{p}}^{\text{s}}}{{\mathbb{C}}}_{abuv}  \pd{ \hat{\xi}_u^{cd}   }{\bar{Y}_v} \dif \bar{\vY}\right).
\eeq
Replacing indices $a$, $b$, $c$ and $d$ with $p$, $q$, $s$ and $t$, Eq~\eqref{Eq:Conformal_cell_stiffness}~is proven, where Eq.~\eqref{Eq:Represent_stiffness_equivalent} has been employed.

\refstepcounter{appendixcounter}
\section*{Appendix \theappendixcounter}
\label{Appendix:Rectangular_solution_Cauchy_Riemann}
The solution of Eq.~\eqref{Eq:BC_determine_solution_rectangle}~are analytic and yield Eq.~\eqref{Eq:Rectangle_solution} is demonstrated in this appendix. Based on the principle of linear superposition, Eq.~\eqref{Eq:Rectangle_lnlmd} is decomposed into three equations as
\beq
\label{Eq:Rectangle_tilde_lnlmd}
\left\{
\begin{aligned}
&\nabla^2\left(\ln\tilde\lambda\right)=0,    \quad  \left(0 \le x_1 \le L, \text{ } 0 \le x_2 \le H\right),\\
&(\ln\tilde\lambda)|_{x_1=0}=\tilde\phi_0\left(x_2\right), \quad(\ln\tilde\lambda)|_{x_1=L}=\tilde\phi_1\left(x_2\right), \\
&(\ln\tilde\lambda)|_{x_2=0}=\tilde\psi_0\left(x_1\right), \quad(\ln\tilde\lambda)|_{x_2=H}=\tilde\psi_1\left(x_1\right), \\
\end{aligned}
\right.
\eeq
\beq
\label{Eq:Rectangle_hat_lnlmd_1}
\;\;\:
\left\{
\begin{aligned}
&\nabla^2\left(\ln\hat\lambda_1\right)=0,    \quad  \left(0 \le x_1 \le L, \text{ } 0 \le x_2 \le H\right),\\
&(\ln\hat\lambda_1)|_{x_1=0}=0,\quad(\ln\hat\lambda_1)|_{x_1=L}=0, \\
&(\ln\hat\lambda_1)|_{x_2=0}=\hat\psi_0\left(x_1\right),\quad(\ln\hat\lambda_1)|_{x_2=H} =\hat\psi_1\left(x_1\right), \\
\end{aligned}
\right.
\eeq
and
\beq
\label{Eq:Rectangle_hat_lnlmd_2}
\;\;\:
\left\{
\begin{aligned}
&\nabla^2\left(\ln\hat\lambda_2\right)=0,    \quad  \left(0 \le x_1 \le L, \text{ } 0 \le x_2 \le H\right),\\
&(\ln\hat\lambda_2)|_{x_1=0}=\hat\phi_0\left(x_2\right),\quad(\ln\hat\lambda_2)|_{x_1=L} =\hat\phi_1\left(x_2\right), \\
&(\ln\hat\lambda_2)|_{x_2=0}=0,\quad(\ln\hat\lambda_2)|_{x_2=H}=0 \\
\end{aligned}
\right.
\eeq
with
\beq
\label{Eq:Rectangle_lnlmd_decompost}
\ln\lambda = \ln\tilde\lambda + \ln\hat\lambda_1 + \ln\hat\lambda_2,
\eeq
and
\beq
\label{Eq:Rectangle_boundary_condition_decomposition}
\left\{
\begin{aligned}
&\phi_0\left(x_2\right)=\tilde\phi_0\left(x_2\right)+\hat\phi_0\left(x_2\right), \quad \phi_1\left(x_2\right)=\tilde\phi_1\left(x_2\right)+\hat\phi_1\left(x_2\right), \\
&\psi_0\left(x_1\right)=\tilde\psi_0\left(x_1\right)+\hat\phi_0\left(x_1\right), \quad \psi_1\left(x_1\right)=\tilde\psi_1\left(x_1\right)+\hat\phi_1\left(x_1\right),
\end{aligned}
\right.
\eeq
where $\tilde\psi_0\left(x_1\right), \tilde\psi_1\left(x_1\right), \tilde\phi_0\left(x_2\right)$ and $\tilde\phi_1\left(x_2\right)$ are all linear functions defined on the boundary of design domain and yield
\beq
\label{Eq:Rectangle_boundary_linear_function}
\left\{
\begin{aligned}
&\tilde\phi_0\left(x_2\right)=\frac{1}{H}\left(\phi_0\left(H\right)- \phi_0\left(0\right)\right)x_2+\phi_0\left(0\right),\\
&\tilde\phi_1\left(x_2\right)=\frac{1}{H}\left(\phi_1\left(H\right)- \phi_1\left(0\right)\right)x_2+\phi_1\left(0\right),\\
&\tilde\psi_0\left(x_1\right)=\frac{1}{L}\left(\psi_0\left(L\right)- \psi_0\left(0\right)\right)x_1+\psi_0\left(0\right),\\
&\tilde\psi_1\left(x_1\right)=\frac{1}{L}\left(\psi_1\left(L\right)- \psi_1\left(0\right)\right)x_1+\psi_1\left(0\right).
\end{aligned}
\right.
\eeq

Firstly, Eq.~\eqref{Eq:Rectangle_tilde_lnlmd} is considered. Since its boundary conditions are all linear, we assume the possible solution of $\ln\lambda$ satisfies
\beq
\label{Eq:Rectangle_Initial_Gauss}
\ln\tilde\lambda = r_1x_1x_2+ r_2x_1+r_3x_2+r_4,
\eeq
where $r_1, r_2, r_3$ and $r_4$ are undetermined coefficients. Substituting Eq.~\eqref{Eq:Rectangle_Initial_Gauss}~into Eq.~\eqref{Eq:Rectangle_tilde_lnlmd} gives
\beq
\label{Eq:Rectangle_tilde_lambda_coefficient}
\begin{aligned}
&r_1 = \frac1{HL}\left(\psi_0\left(0\right)-\psi_0\left(L\right)+\psi_1\left(L\right)-\psi_1\left(0\right)\right),&\quad &r_2=\frac1{L} \left( \psi_0\left(L\right)-\psi_0\left(0\right) \right),\\
&r_3 = \frac1{H}\left( \psi_1\left(0\right)-\psi_0\left(0\right) \right),&\quad &r_4 = \psi_0\left(0\right).
\end{aligned}
\eeq

Then Solving Eqs.~\eqref{Eq:Rectangle_hat_lnlmd_1}~and~\eqref{Eq:Rectangle_hat_lnlmd_2}~directly by using the method of separation of variables, it is obtained that
\begin{gather}
\ln\hat\lambda_1 = \sum_{k=1}^\infty \left(a_k\cosh\frac{k\pi{x_2}}{L}+b_k\sinh\frac{k\pi{x_2}}{L}\right)\sin\frac{k\pi{x_1}}{L}
,\label{Eq:Rectangle_hat_lnlmd_solution_1}\\
\ln\hat\lambda_2 = \sum_{k=1}^\infty \left(c_k\cosh\frac{k\pi{x_1}}{H}+d_k\sinh\frac{k\pi{x_1}}{H}\right)\sin\frac{k\pi{x_2}}{H}, \label{Eq:Rectangle_hat_lnlmd_solution_2}
\end{gather}
where $a_k,b_k,c_k$ and $d_k$ are arbitrary constants and yield
\beq
\label{Eq:Rectangle_boundary_fourier_series_coefficients}
\begin{aligned}
&a_k=\frac{2}{L}\int_0^L\hat\psi_0\left({x_1}\right)\sin\frac{k\pi{x_1}}{L}\dif{x_1},  \\
&b_k=\frac{2}{L\sinh\frac{k\pi{H}}{L}}\int_0^L\left(\hat\psi_1\left(x_1\right)-\hat\psi_0({x_1) \cosh\frac{k\pi{H}}{L}}\right)\sin\frac{k\pi{x_1}}{L}\dif{x_1},\\
&c_k=\frac{2}{H}\int_0^H\hat\phi_0\left({x_2}\right)\sin\frac{k\pi{x_2}}{H}\dif{x_2},\\
&d_k=\frac{2}{H\sinh\frac{k\pi{L}}{H}}\int_0^L\left(\hat\phi_1\left(x_2\right)-\hat\phi_0({x_2) \cosh\frac{k\pi{L}}{H}}\right)\sin\frac{k\pi{x_2}}{H}\dif{x_2},
\end{aligned}
\eeq
respectively.
Incorporating Eqs.~\eqref{Eq:Rectangle_Initial_Gauss}--\eqref{Eq:Rectangle_boundary_fourier_series_coefficients} into Eq.~\eqref{Eq:Rectangle_lnlmd_decompost}~gives the expression for $\ln\lambda$ in Eq.\eqref{Eq:Rectangle_lnlmd_solution}. Similarly, we can find Eq.~\eqref{Eq:Rectangle_theta_solution} is the solution of Eq.~\eqref{Eq:Rectangle_theta}.

\end{document}